\documentclass[12pt]{article}
 \usepackage{amsmath, amsthm, amssymb}
 \usepackage[margin=1in]{geometry}

\newtheorem{theorem}{Theorem}
\newtheorem{lemma}[theorem]{Lemma}

\newif\ifpdf
	\ifx\pdfoutput\undefined
	\pdffalse 
	\else
	\pdfoutput=1 
	\pdftrue
	\fi

	\ifpdf
	\usepackage[pdftex]{graphicx}
	\else
	\usepackage{graphicx}
	\fi

\begin{document}

\ifpdf
	\DeclareGraphicsExtensions{.pdf, .png, .jpg, .tif}
	\else
	\DeclareGraphicsExtensions{.eps, .jpg}
	\fi

\begin{center}
{\Large {\bf Bayes variable selection in semiparametric linear models\\}} 
Suprateek Kundu$^1$ and David B. Dunson$^2$ \\
\end{center}

{\noindent {\bf \quad Abstract:}}	There is a rich literature proposing methods and establishing asymptotic properties of Bayesian variable selection methods for parametric models, with a particular focus on the normal linear regression model and an
increasing emphasis on settings in which the number of candidate predictors ($p$) diverges with sample size ($n$).  Our focus is on generalizing methods and asymptotic theory established for mixtures of $g$-priors to semiparametric linear regression models having unknown residual densities. Using a Dirichlet process location mixture for the residual density, we propose a
semiparametric $g$-prior which incorporates an unknown matrix of cluster allocation indicators.  For this class of priors, posterior computation can proceed via a straightforward stochastic search variable selection algorithm.  In addition, Bayes factor and variable selection consistency is shown to result under various cases including proper and improper priors on $g$ and $p>n$, with the models under comparison restricted to have model dimensions diverging at a rate less than $n$. 

	
\begin{footnotesize}
$^1$Suprateek Kundu is a doctoral candidate in the Dept. of Biostatistics, UNC Chapel Hill, Chapel Hill, NC 27599, USA (skundu@email.unc.edu).\\
$^2$David B. Dunson is professor in Dept. Statistical Science, Duke University, Durham, NC 27708, USA (dunson@stat.duke.edu).\\
\end{footnotesize}

\newpage
	
{\noindent{\bf 1.\quad INTRODUCTION}  }  \\
Bayesian variable selection is very widely applied, with a rich literature on alternative priors and computational methods.  For a recent review of Bayesian variable selection methods, refer to O'Hara and Sillanp\"a\"a (2009).  Most of the literature has focused on Gaussian linear regression models, with common methods including stochastic search variable selection (SSVS) (George and McCulloch, 1993; 1997), reversible jump MCMC (Green, 1995) and adaptive shrinkage (Tibshirani, 1996; Park and Casella, 2008; Yi and Xu, 2008).  Such methods can be applied directly for kernel or basis function selection in nonlinear regression with Gaussian residuals (Smith and Kohn, 1996) and can be adapted to accommodate generalized linear models with outcomes in the exponential family (Raftery and Richardson 1993; Meyer and Laud 2002).  

It is well known that Bayesian variable selection can be sensitive to the prior, and there is an increasingly rich literature showing asymptotic properties providing support for carefully-chosen priors, such as mixtures of g-priors (Zellner and Siow, 1980; Liang et. al., 2008), with such priors also having appealing computational properties.  This literature is essentially entirely focused on Gaussian linear regression models, and the emphasis of this article is on developing methods that generalize this work to semiparametric regression models having unknown residual distributions.

To set the stage, first consider the well-studied problem of comparison of linear models of the following type: 
\begin{eqnarray}
M_1: Y^n &=& \alpha 1_n + X_{\gamma_1}\beta_{\gamma_1} + \epsilon_1, \quad \epsilon_1 \sim N(0,\tau^{-1} I_n), \nonumber \\
M_2: Y^n &=& \alpha 1_n + X_{\gamma_2}\beta_{\gamma_2} + \epsilon_2, \quad \epsilon_2 \sim N(0,\tau^{-1} I_n), \label{eq:LRM}
\end{eqnarray}
where Y$^n$ is n$\times$1 vector of responses, $\alpha$ is the common intercept, X$_{\gamma_j}$ is a n$\times p_j$ design matrix (j=1,2) excluding the column of intercepts, and $\epsilon_j$'s are Gaussian residuals, j=1,2. The models may or may not be nested, and the number of candidate predictors is $p$. Among numerous model selection criteria available for such comparisons, the Bayes factor (Kass and Raftery, 1995) has received substantial attention as the most widely accepted Bayesian measure of the weight of evidence in the data in favor of one model over another. The Bayes factor for comparing $M_1$ versus $M_2$ based on a sample Y$^n$ is defined as BF$^n_{12}=\frac{L(Y^n|M_1)}{L(Y^n|M_2)}$, the ratio of marginal likelihoods under $M_1$ and $M_2$. Assuming one of the models under comparison is true, Bayes factor consistency refers to the phenomenon where BF$^n_{12}\stackrel{P}{\rightarrow} \infty$ as $n \to \infty$ under $M_1$ and BF$^n_{12}\stackrel{P}{\rightarrow} 0$ as $n \to \infty$ under $M_2$. A stronger form of consistency is also possible when the convergence happens almost surely. When comparing the true model pairwise to each model in a list, Bayes factor consistency typically implies that the posterior probability on the true model goes to one.

Although priors most commonly used in practice assume {\em a priori} independence in the elements of the coefficient vectors ($\beta_1$ and $\beta_2$), priors that have been shown to result in Bayes factor consistency typically incorporate dependence.  Examples include the intrinsic prior (Berger and Pericchi, 1996; Moreno, Bertolino and Racugno, 1998), and Zellner's $g$-prior (Zellner, 1986) specified by $\beta_j \sim N(0,g\tau^{-1}(X_j'X_j)^{-1})$, j=1,2. The intrinsic  priors have proven to behave very well for multiple testing problems (Casella and Moreno, 2006). Zellner's $g$-prior allows for a convenient correlation structure and can control for the amount of prior information relative to the sample through only one hyperparameter $g$. Among others, Fern$\acute{a}$ndez et al. (2001) investigated Bayes factor consistency under various choices of fixed $g$, which was allowed to depend on the sample size and/or the number of candidate predictors. In order to resolve difficulties associated with a fixed choice of g, such as Bartlett's paradox (Bartlett, 1957; Jeffreys 1961) and information paradox (Zellner 1986; Berger and Pericchi 2001), Zellner and Siow (1980) placed an inverse-gamma prior on $g$, while Liang et. al. (2008) extended the idea of Strawderman (1971) to the regression context by proposing hyper-$g$ and hyper-$g/n$ priors on $g$, under which they established Bayes factor consistency. The above approaches entail specifying improper priors on common model parameters and proper priors on model parameters unique to any one model, which results in a prior specification for the more complex model depending upon the simpler model. To avoid such pitfalls, Guo and Speckman (2009) adopted the idea of Marin and Robert (2007) and placed mixtures of $g-$priors on all the elements of both $\beta_1$ and $\beta_2$, which leads to tractable Bayes factors as well as Bayes factor consistency. 

There has also been a growing interest in model selection procedures for normal linear models when the number of candidate predictors ($p$) increase with sample size ($n$). Such increases occur in a wide variety of applications, such as in nonparametric regression when the number of candidate kernels or basis functions depends on $n$. Shao (1997) analyzed the consistency of several frequentist and Bayesian approximation criteria for model selection in normal linear models with increasing model dimensions, assuming the true model to be the submodel minimizing the average squared prediction error. Moreno et. al. (2010) examined consistency of Bayes factors and the BIC under intrinsic priors for nested normal linear models, when the dimension of the parameter space increases with the sample size. Jiang (2007) considered Bayesian variable selection in generalized linear models in $p>n$ settings and provided conditions to obtain near optimal rates of convergence in estimating the conditional predictive distribution, but did not consider asymptotic properties in selecting the important predictors.

To our knowledge, this area has entirely focused on parametric models with a particular focus on normal linear regression. Such a parametric assumption on the residual error is rather stringent and may not hold in practice, thus invalidating the earlier assumption of the true model belonging to the class of models under comparison and potentially leading to inconsistent Bayes factors. In Section 5, simulations illustrate that when residuals are generated from a bimodal distribution, Bayesian variable selection under a Gaussian linear regression model tends to have poor performance.  With this motivation, our focus is on developing Bayes variable selection methods that do not require Gaussian residuals and that can be shown to be consistent. 

There is a limited literature on variable selection in Bayesian regression models having unknown residual distributions.  Kuo and Mallick (1997) consider an accelerated failure time model for time-to-event data containing a linear regression component and a mixture of Dirichlet processes for the residual density.  To perform variable selection, they add indicator variables to the regression function and implement an MCMC algorithm.  Also, in the survival analysis setting, Dunson and Herring (2005) proposed a Bayesian approach for selecting predictors in a semiparametric hazards model that allows uncertainty in whether predictors enter in a multiplicative or additive manner.  Kim, Tadesse and Vannucci (2006) instead define a Bayesian variable selection approach, which uses a Dirichlet process to define clusters in the data, while updating the variable inclusion indicators using a Metropolis scheme.  Mostofi and Behboodian (2007) models a symmetric and unimodal residual density using a Dirichlet process scale mixtures of uniforms, while conducting Bayesian variable selection.  Chung and Dunson (2009) modeled the conditional response density given predictors using a flexible probit stick-breaking mixture of Gaussian linear models, allowing variable selection via a Bayesian stochastic search method.  

These articles focused on defining methodology and computational algorithms, but without study of theoretical properties, such as consistency.  In fact, to our knowledge, there has been no previous work on consistent Bayesian variable selection in semi-parametric models, though there is recent work on consistent non-parametric Bayesian model selection (Ghosal, Lember and van der Vaart, 2008 among others). It is not straightforward to apply such theory directly to the problem of variable selection in semiparametric linear models.  

With this motivation, we define a practical, useful and general methodology for Bayesian variable selection in semiparametric linear models, while providing basic theoretical support by showing Bayes factor and variable selection consistency. We accomplish this by generalizing the methods and asymptotic theory for mixtures of $g$-priors to linear regression models with unknown residuals characterized via Dirichlet process (DP) location mixture of Gaussians. We propose a new class of mixtures of semi-parametric $g$-priors, which results in consistent Bayesian variable selection even when there are many more candidate predictors ($p$) than samples ($n$) as long as the prior assigns probability zero to models having greater than or equal to $n$ predictors.  Additionally, posterior computation for the proposed method is straightforward via an SSVS algorithm. 

Section 2 develops the proposed framework. Section 3 considers asymptotic properties. Section 4 outlines algorithms for posterior computation.  Section 5 contains a simulation study. Section 6 applies the approach to a type 2 diabetes data example, and Section 7 discusses the results.

\vskip 12pt

{\noindent {\bf 2.\quad MIXTURES OF SEMIPARAMETRIC $g$-PRIORS} }\\
{\noindent {\bf 2.1 \quad MODEL FORMULATION} }\\
In this section, we propose a new class of priors for Bayesian variable selection in linear regression models with an unknown residual density characterized via a Dirichlet process (DP) location mixture of Gaussians. In particular, let 
\begin{eqnarray}
y_i &=& { x}_{\gamma,i}'\beta_{\gamma} + \epsilon_i,\quad \epsilon_i \sim f, i=1,\ldots,n, \nonumber \\
f(\cdot) &=& \int N(\cdot;\alpha,\tau^{-1}) dP(\alpha), \quad P \sim DP(mP_0), \quad P_0 = N(0,\tau^{-1}), \label{eq:base} 
\end{eqnarray}
where $\gamma = (\gamma^1,\ldots, \gamma^p)' \in \Gamma$ is a vector of variable inclusion indicators, with $\gamma^j$=I($j$th predictor is included in the model) and $\sum_{j=1}^p\gamma^j=p_\gamma$, $\beta_{\gamma} = \{ \beta_j:\gamma^j = 1, j=1,\ldots, p$\}, ${ x}_{\gamma,i} = \{ x_{ij}: \gamma^j = 1, j=1,\ldots, p\} \in \mathcal{X}$ and does not include an intercept, and $f$ is a density with respect to Lebesgue measure on $\Re$. For simplicity, we model $f$ as having an unknown mean instead of including an intercept $\alpha$ as in (\ref{eq:LRM}). The number of candidate predictors $p$ may or may not increase with the sample size $n$. We can address the prior uncertainty in subset selection by placing a prior on $\gamma$, while the prior on $\beta_{\gamma}$ characterizes prior knowledge of the size of the coefficients for the selected predictors. 

The DP mixture prior on the density $f$ induces clustering of the $n$ subjects into $k$ groups, with each group having a distinct intercept in the linear regression model.  Let A denote an $n \times k$ allocation matrix, with A$_{ij}$ = 1 if the $i$th subject is allocated to the $j$th cluster and 0 otherwise.  The $j$th column of $A$ then sums to $n_j$, the number of subjects allocated to cluster $j$, with $\sum_{j=1}^k n_j = n$. Following Kyung, Gill and Casella (2009), conditionally on the allocation matrix A, (\ref{eq:base}) can be represented as the linear model 
\begin{eqnarray}
Y^n = A\eta + X_\gamma \beta_\gamma + \epsilon, \quad \eta \sim N(0,\tau^{-1} I_k), \quad \epsilon \sim N(0,\tau^{-1} I_n), \label{eq:matbase}
\end{eqnarray}
where $X_{\gamma} = ( x_{\gamma,i}, i=1,\ldots, n )'$. 

In keeping with the mixtures of $g$-priors literature, we would like the prior on the regression coefficients to retain the essential elements of Zellner's $g$-prior, but at the same time to be suitably adapted to reflect the semi-parametric nature of the model in question - more specifically, the clustering of responses by the DP kernel mixture prior. To this effect, we propose a mixture of semi-parametric $g$-priors which is constructed to scale the covariance matrix in Zellner's $g$-prior to reflect the clustering phenomenon as follows:
\begin{eqnarray}
\pi(\beta_{\gamma})= N( 0, g\tau^{-1}(X_{\gamma}' \Sigma_A^{-1} X_{\gamma})^{-1} ), \quad \Sigma_A = I+AA', \quad g \sim \pi(g). \label{eq:beta}
\end{eqnarray}
Prior (\ref{eq:beta}) inherits the advantages of the traditional mixtures of $g$-priors including computational efficiency in computing marginal likelihoods (conditional on A) and robustness to mis-specification of $g$.  In addition, the prior can be interpreted as having arisen from the analysis of a conceptual sample generated using a scaled design matrix $\Sigma_A^{-1/2}X_{\gamma}$, reflecting the clustering phenomenon due to the DP kernel mixture prior. Moreover, the proposed prior leads to Bayes factor and variable selection consistency in semi-parametric linear models (\ref{eq:base}), as highlighted in the sequel. 

Note that since $(X_\gamma' \Sigma_A^{-1} X_\gamma)^{-1}\ge(X_\gamma'X_\gamma)^{-1}$ for any allocation matrix A, the prior variance of Y conditional on (A,$g,\tau^{-1}$) is higher for the semi-parametric $g-$prior as compared to the traditional $g-$prior. To assess the influence of A on the prior for $\beta_\gamma$, we did simulations which revealed that for fixed ($n$, $p$), var($\beta_k$) increases but the cov($\beta_k,\beta_l$) decreases as the number of underlying clusters in the data increase ($k,l=1,\ldots,p, k\ne l$). This suggests that as the number of clusters increase, the components of $\beta$ are likely to be more dispersed with decreasing association between each other. 

\vskip 12pt

{\noindent {\bf 2.2 \quad Bayes Factor in Semiparametric Linear Models}} \\
Throughout the rest of the paper, we will assume that the data Y$^n$=$(Y_1,\ldots,Y_n)'$ are generated from the true model $\mathcal{M}_T: Y^n = X_{\gamma_1}\beta_{\gamma_1} + \epsilon,$ with $\epsilon_{i}$ i.i.d. from the true residual density $f_0$, which is a density on $\Re$ with respect to Lesbesgue measure. For modeling purposes, we put a DP location mixture of Gaussians prior on the unknown $f_0$. For pairwise comparison, we evaluate the evidence in favor of $\mathcal{M}_1$ compared to $\mathcal{M}_2$ using Bayes factor, where 
\begin{eqnarray}
\mathcal{M}_1 &:& Y^n = X_{\gamma_1} \beta_{\gamma_1}  + \epsilon_1, \quad \epsilon_{1i} \sim f \nonumber \\
\mathcal{M}_2 &:& Y^n = X_{\gamma_2} \beta_{\gamma_2}  + \epsilon_2, \quad \epsilon_{2i} \sim f \nonumber \\
f(\cdot) &=& \int N(\cdot;\alpha,\tau^{-1}) dP(\alpha), \quad P \sim DP(mP_0), \quad P_0 = N(0,\tau^{-1}) \nonumber \\
\beta_{\gamma_j} &\sim& \pi(\beta_{\gamma_j}), j=1,2, \quad \pi(\tau^{-1}) \propto 1/\tau^{-1}, \quad g \sim \pi(g),  \label{eq:SPLM}
\end{eqnarray}
where $\gamma_j$ indexes models of dimension $p_j$ in the model space $\mathcal{M}$ (j=1,2) and $\pi(\beta_{\gamma_j})$ is defined in (\ref{eq:beta}). Our prior specification philosophy is similar to the one adopted by Guo and Speckman (2009) for normal linear models, in that we assign proper priors on all elements of both $\beta_{\gamma_1},\beta_{\gamma_2}$ conditional on $(g,\tau^{-1})$, and an improper prior on $\tau^{-1}$ (for a more objective assessment). However unlike Guo and Speckman (2009), our focus is on Bayesian variable selection in semi-parametric linear models. 

Note that the likelihood of the response after marginalizing out $\eta$ in (\ref{eq:matbase}) turns out to be L$(Y^n|A,\beta_\gamma,\tau^{-1})= N(X_\gamma\beta_\gamma,\tau^{-1} \Sigma_A)$ (Kyung et. al., 2009). Thus conditional on A, $Z_A = \Sigma_A^{-1/2}Y^n \sim N(\Sigma_A^{-1/2}X_\gamma\beta_\gamma,\tau^{-1} I_n)$, and we are in the normal linear model set-up: 
\begin{eqnarray}
Z_A = \tilde{X}_{A,\gamma}\beta_\gamma + \epsilon, \quad \epsilon \sim N(0, \tau^{-1} I_n), \quad \pi(\beta_{\gamma})= N( 0, g\tau^{-1}(\tilde{X}'_{A,\gamma} \tilde{X}_{A,\gamma})^{-1} ),
\end{eqnarray}
where  $\tilde{X}_{A,\gamma}=\Sigma_A^{-1/2}X_\gamma$. Under a mixture of semi-parametric $g$-priors, we can directly use expression (17) in Guo and Speckman (2009) to obtain for j=1,2
\begin{eqnarray}
L(Z_A|\mathcal{M}_j)\equiv L(Y^n|A,\mathcal{M}_j) \propto (Z_A'Z_A)^{-n/2} \int_0^\infty (1+g)^{-p_j/2}\left[ 1-\frac{g}{1+g} \frac{Z_A'\tilde{H}_{A,j}Z_A}{Z_A'Z_A}\right]^{-n/2} \pi(dg), \label{eq:marginal}
\end{eqnarray}
where $\tilde{H}_{A,j}=\tilde{X}_{A,\gamma_j}(\tilde{X}_{A,\gamma_j}'\tilde{X}_{A,\gamma_j})^{-1}\tilde{X}_{A,\gamma_j}'$, the equivalent of a hat matrix in standard linear regression. Also, marginalizing over all possible subcluster allocations for a given sample size n, the following form for the marginal likelihood can be obtained (Kyung et. al., 2009):
\begin{eqnarray}
L(Y^n|\mathcal{M}_j) = \frac{\Gamma (m)}{\Gamma (m+n)}\sum_{k=1}^n m^k \sum_{A \in \mathcal{A}_k} \prod_{i=1}^k \Gamma(n_i) L(Y^n|A,\mathcal{M}_j) = \sum_{A_l\in \mathcal{C}_n} w_l L(Z_{A_l}|\mathcal{M}_j), \label{eq:kyung}
\end{eqnarray}
where $\mathcal{A}_k$ is the collection of all possible n$\times$k matrices corresponding to different allocations of n subjects into k subclusters, $\mathcal{C}_n$ is the collection of all possible allocation matrices for a sample size n with $\sum_{A_l\in \mathcal{C}_n} w_l=1$. In the limiting case as $n\to \infty$, we have $\mathcal{C}_\infty$ as the class of `limiting allocation matrices'. Using (\ref{eq:marginal}), the Bayes factor in favor of $\mathcal{M}_2$ conditional on the allocation matrix A is given by
\begin{eqnarray}
BF^n_{21,A}&=& \frac{L(Z_{A}|\mathcal{M}_2)}{L(Z_{A}|\mathcal{M}_1)} = \frac{\int_0^\infty (1+g)^{-p_2/2}\left[1-\frac{g}{1+g}\tilde{R}_{A,2}^2\right]^{-n/2}\pi(dg)}{\int_0^\infty (1+g)^{-p_1/2}\left[1-\frac{g}{1+g}\tilde{R}_{A,1}^2\right]^{-n/2}\pi(dg)} , \label{eq:BFA}
\end{eqnarray}
where $\tilde{R}_{A,j}^2= Z_A'\tilde{H}_{A,j}Z_A/Z_A'Z_A$, (j=1,2). Finally using (\ref{eq:kyung}), the unconditional Bayes factor marginalizing out A in favor of $\mathcal{M}_2$ is 
\begin{eqnarray}
BF^n_{21} &=& \frac{L(Y^n|\mathcal{M}_2)}{L(Y^n|\mathcal{M}_1)}= \frac{ \sum_{A_l\in \mathcal{C}_n} w_l L(Z_{A_l}|\mathcal{M}_2)}{\sum_{A_l\in \mathcal{C}_n} w_l L(Z_{A_l}|\mathcal{M}_1)}. \label{eq:BF}
\end{eqnarray}

\vskip 12pt

{\noindent{\bf 3. \quad ASYMPTOTIC PROPERTIES}}\\
In this section we focus on asymptotic properties including Bayes factor and variable selection consistency. Before proceeding, note that the standard assumptions made for establishing Bayes factor consistency in linear models (\ref{eq:LRM}) are: \\
\emph{(A1$'$)} $\lim_{n \to \infty}\frac{\beta_{\gamma_1}'(X_{\gamma_1}'X_{\gamma_1})\beta_{\gamma_1}}{n} \to b_{1}>0$ under $M_1$, \\
\emph{(A2$'$)} If $M_1\not\subseteq M_2$, $0 \le \lim_{n\to \infty}\frac{\beta_{\gamma_1}'X_{\gamma_1}'H_{2}X_{\gamma_1}\beta_{\gamma_1}}{n} \to b_{2}$, $0 \le b_2<b_{1}$ under $M_1$,\\
where $H_2=X_{\gamma_2}(X_{\gamma_2}'X_{\gamma_2})^{-1}X_{\gamma_2}'$, the hat matrix in $M_2$. A necessary condition for assumption \emph{(A1$'$)} to hold is that X$_{\gamma_1}$  has full rank, which is likely to be satisfied for fixed model dimensions but can not be guaranteed for increasing model dimensions without further assumptions. Conditional on the limiting allocation matrix A$\in \mathcal{C}_\infty$, we make similar assumptions. For fixed $p_j$ and conditional on A$\in \mathcal{C}_\infty$, we assume \\
\emph{(A1)} $\lim_{n \to \infty}\frac{\beta_{\gamma_1}'(X_{\gamma_1}'\Sigma_A^{-1}X_{\gamma_1})\beta_{\gamma_1}}{n} \to b_{A,1}>0$ under $\mathcal{M}_1$. \\
\emph{(A2)}: If $\mathcal{M}_1 \not\subseteq \mathcal{M}_2$, $\lim_{n\to \infty}\frac{\beta_{\gamma_1}'\tilde{X}_{A,\gamma_1}'\tilde{H}_{A,2}\tilde{X}_{A,\gamma_1}\beta_{\gamma_1}}{n} \to b_{A,2}\in [0,b_{A,1})$ under $\mathcal{M}_1$.  \\
For $p_j=O(n^{a_j})$ (j=1,2) with $0\le a_1<a_2<1$, conditional on A$\in \mathcal{C}_\infty$ we assume \emph{(A1)} and \\
\emph{($\tilde{A}$2)}: If $\mathcal{M}_1 \not\subseteq \mathcal{M}_2$, $\lim_{n\to \infty}\frac{\beta_{\gamma_1}'\tilde{X}_{A,\gamma_1}'\tilde{H}_{A,2}\tilde{X}_{A,\gamma_1}\beta_{\gamma_1}}{n} \to b_{A,2}\in (0,b_{A,1})$ under $\mathcal{M}_1$.  \\
Note that \emph{(A1)}$\Rightarrow$\emph{(A1$'$)} (as $X_\gamma' \Sigma_A^{-1} X_\gamma\le X_\gamma'X_\gamma$), so that assumption \emph{(A1)} is a stronger version of \emph{(A1$'$)}. Further, for the two extreme cases when A=$1_n$ and A=I$_n$, \emph{(A1$'$)}$\Rightarrow$\emph{(A1)}. To see this, note that $X_{\gamma_1}'\Sigma_{A=1_n}^{-1}X_{\gamma_1} \approx X_{\gamma_1}'X_{\gamma_1} - n\bar{X}_{\gamma_1}'\bar{X}_{\gamma_1}$, where $\bar{X}_{\gamma_1}$ is a $1 \times p$ vector containing the column means of X$_{\gamma_1}$. This implies $\frac{\beta_{\gamma_1}'(X_{\gamma_1}'\Sigma_{A=1_n}^{-1}X_{\gamma_1})\beta_{\gamma_1}}{n}\approx \frac{\beta_{\gamma_1}'(X^{c'}_{\gamma_1}X^c_{\gamma_1})\beta_{\gamma_1}}{n}>0$ (plugging in X$^c_\gamma$ for X$_\gamma$ in \emph{(A1$'$)}), where X$^c_{\gamma_1}$ is the centered version of $X_{\gamma_1}$ such that $1'_nX^c_{\gamma_1}=0_{1\times p}$. On the other hand for A=I$_n$, we have $\frac{X_{\gamma_1}'\Sigma_A^{-1}X_{\gamma_1}}{n} = \frac{1}{2}\frac{X_{\gamma_1}'X_{\gamma_1}}{n}>0$. Assumptions \emph{(A2)}, \emph{($\tilde{A}$2)} can be interpreted as a positive `limiting distance' between two models corresponding to design matrices X$_{\gamma_1}$ and X$_{\gamma_2}$  in (\ref{eq:matbase}) conditional on A$\in \mathcal{C}_\infty$, after marginalizing out $\eta$, i.e. $\Delta_{21,A}=\lim_{n\to \infty}\frac{\beta_{\gamma_1}'\tilde{X}_{A,\gamma_1}'(I_n-\tilde{H}_{A,2})\tilde{X}_{A,\gamma_1}\beta_{\gamma_1}}{n\tau^{-1}}=\frac{b_{A,1}-b_{A,2}}{\tau^{-1}} \in (0,\infty)$. Such a `limiting distance' ($\Delta_{21,A}$) can be considered as a natural extension of the definition of distance between two normal linear models in Casella et. al. (2009) and Moreno et. al. (2010).   

The following lemma gives the limits of quantities such as $\tilde{R}_{A,j}^2= Z_A'\tilde{H}_{A,j}Z_A/Z_A'Z_A$ (A$\in \mathcal{C}_\infty$, j=1,2), which would be useful for establishing asymptotic properties. The proof follows directly from the fact that conditional on allocation matrix A, $Z_A = \Sigma_A^{-1/2}Y^n \sim N( \tilde{X}_{A,\gamma}\beta, \tau^{-1} I_n)$, and using Lemmas 1 and 2 of Guo and Speckman (2009). From here on, we shall make all probability statements under the model $\mathcal{M}_1$ as defined in (\ref{eq:SPLM}).

\begin{lemma} 
Suppose assumptions (A1), (A2) and ($\tilde{A}$2) hold. \\
(i) If $\mathcal{M}_1\subset\mathcal{M}_2$, then conditional on A$\in \mathcal{C}_\infty$, $\tilde{R}^2_{A,1}\stackrel{a.s.}{\rightarrow}\frac{b_{A,1}}{\tau^{-1}+b_{A,1}}, \quad \tilde{R}^2_{A,2}\stackrel{a.s.}{\rightarrow}\frac{b_{A,1}}{\tau^{-1}+b_{A,1}}$. \\
(ii) If $\mathcal{M}_1\not\subseteq\mathcal{M}_2$, then conditional on A$\in \mathcal{C}_\infty$, $\tilde{R}^2_{A,1}\stackrel{a.s.}{\rightarrow}\frac{b_{A,1}}{\tau^{-1}+b_{A,1}}, \quad \tilde{R}^2_{A,2}\stackrel{a.s.}{\rightarrow}\frac{b_{A,2}}{\tau^{-1}+b_{A,1}}$ .
\end{lemma}

Although the next result establishes Bayes factor consistency in semi-parametric linear models (\ref{eq:SPLM}) under the class of proper priors for $g$, the result can be extended to improper priors $\pi(g) \propto \frac{1}{1+g}$. For fixed $p$, $p_1$ can be greater or less than $p_2$, while for increasing $p$ we compare models with $p_1=O(n^{a_1})$, $p_2=O(n^{a_2})$ and $0 \le a_1<a_2<1$, which involves the special case of fixed $p_1$ but increasing $p_2$. As elaborated in Guo and Speckman (2009), the class of priors $\pi(g)$ considered here include hyper-$g$ ($\frac{a-2}{2}(1+g)^{-a/2}$) and hyper-$g/n$ ($\frac{a-2}{2n}(1+g/n)^{-a/2}$) priors, with $2<a\le 4$ (Liang et. al. 2008), Zellner-Siow and beta-prime priors. Let the notation $a_n\approx b_n$ imply that $\lim_{n\to \infty} a_n/b_n>0$ almost surely. We assume the following conditions on $\pi(g)$: \\
\emph{(A3):} There exists a constant k$\ge0$ such that $\int_{a_n}^{c_0a_n}\pi(dg)\approx n^{-k}$ for any constant $c_0>1$ and any sequence $a_n\approx n$. \\
\emph{(A4):} There exists a constant k$_u$ such that k-($p_2-p_1)/2< k_u \le$ k and $\int_0^\infty (1+g)^{k_u}\pi(dg) \approx 1$. \\
Assumption \emph{(A3)} ensures that the prior mass for the tail decreases exponentially fast, which is a weak condition and quite reasonable. 
\newtheorem*{thmI}{Theorem I}
\begin{thmI}
Suppose assumptions (A1), (A2) and ($\tilde{A}$2) hold . \\
(I) Suppose $p_1$ and $p_2$ are fixed. If $\mathcal{M}_1\subset\mathcal{M}_2$, then under $\mathcal{M}_1$ and assumptions (A3), (A4), BF$^n_{21}\stackrel{P}{\rightarrow}0$ as n$\to \infty$ and if $p_2-p_1 > 2+2(k-k_u)$, BF$^n_{21}\stackrel{a.s.}{\rightarrow}0$ as n$\to \infty$. Further, if $\mathcal{M}_1\not\subseteq\mathcal{M}_2$, then under $\mathcal{M}_1$ and assumption (A3), BF$^n_{21}\stackrel{a.s.}{\rightarrow} 0$ as n$\to \infty$. \\
(II)  Suppose $p_1=O(n^{a_1})$ and $p_2=O(n^{a_2})$, with $0 \le a_1<a_2<1$. Then under $\mathcal{M}_1$ and assumption (A3), BF$^n_{21}\stackrel{a.s.}{\rightarrow}0$ as n$\to \infty$. 
\end{thmI}
\begin{footnotesize} \textbf{REMARK 1}. \end{footnotesize} The above result can be easily extended to improper priors $\pi(g)\propto \frac{1}{1+g}$.

In settings in which there are not two models under consideration but many, often it is of interest to see if the posterior model probability P($\mathcal{M}_1|Y^n$) goes to 1 as n$\to \infty$. The next theorem gives such a result making use of a sequence of prior model probabilites depending on $n$ and assuming that the growth rate of $\mathcal{M}_1$ is known, for increasing model dimensions.
\newtheorem*{thmII}{Theorem II}
\begin{thmII}
Suppose the conditions of Theorem I hold. For fixed $p$ and under $\mathcal{M}_1$, P($\mathcal{M}_1|Y^n$)$\stackrel{}{\to}$ 1 for any $\left\{\pi(\mathcal{M}_\gamma):\gamma \in \Gamma,\pi(\mathcal{M}_1)>0\right\}$. For increasing $p$ ($\ge p_1$) and under $\mathcal{M}_1$, P($\mathcal{M}_1|Y^n$) $\to$ 1 for $\left\{\pi^n(\gamma_{j_{l}}) \propto  2^{-p_j/2}I[p_j\le O(n^{a_1})] +  N_j^{-1}I[O(n^{a_1})<p_j\le (n-1)\wedge p]\right\}$, where $\gamma_{j_{l}}$ denotes the lth model having $p_j$ predictors, $l=1,\ldots,N_j$, with $N_j = \binom{p}{p_j}$. 
\end{thmII}
\begin{footnotesize} \textbf{REMARK 2}. \end{footnotesize} The mode of convergence of P($\mathcal{M}_1|Y^n$) under $\mathcal{M}_1$ is the same as that of the associated conditional Bayes factors. 

\vskip 12pt

{\noindent{\bf 4. \quad POSTERIOR COMPUTATION}}\\
We propose an MCMC algorithm for posterior computation, which combines a stochastic search variable selection algorithm (George and McCulloch, 1997) with recently proposed methods for efficient computation in Dirichlet process mixture models.  In particular, we utilize the slice sampler of Walker (2007) incorporating the modification of Yau et al. (2011).  Using Sethuraman's (1994) stick-breaking representation, let 
\begin{eqnarray*}
P = \sum_{j=1}^{\infty} w_j \alpha_j, \quad \alpha_j \sim N(0,\tau^{-1}), \quad w_j = \nu_j \prod_{l<j}(1-\nu_l), \quad \nu_l \sim Beta(1,m).
\end{eqnarray*}
The slice sampler of Walker (2007) relies on augmentation with uniform latent variables as follows 
\begin{eqnarray*}
f_{w,\alpha}(y)  = \sum_{j \in B_w(u)} N(y|\alpha_j), \quad B_w(u)= \left\{j:w_j>u\right\}.
\end{eqnarray*}
For sampling the DP precision parameter, we specify the prior $m \sim Ga(a_m,b_m)$, as such a hierarchical specification is likely to ensure better performance by increasing the support of the prior. Further, we assume equal prior inclusion probability for all predictors, i.e. $\pi(\gamma^k=1)= \frac{1}{2}$, k=$1,\ldots,p$. We outline the posterior computation steps briefly below:\\
{\em Step 1.1:} Update the $\nu'$s after marginalizing out the augmented uniform variable using $\pi(\nu_h|-)=Be(1+n_h,\sum_{j>h}n_j + m)$. \\
{\em Step 1.2:} Update the augmented uniform variables from its full conditional as described in Walker (2007). \\
{\em Step 2:} Update the allocation of atoms to different subjects using $f(y_i|u_i,S_i=h)\propto N(y_i|\alpha_h,x_{\gamma,i},\beta_\gamma,\tau^{-1})I(h\in B_w(u_i))$, h=1,\ldots,M\\
{\em Step 3:} Update the precision parameter of the DP using $\pi(m|-)=Ga(a_m + M,b_m - \sum_{l=1}^M\log(1-\nu_l) )$, where M is the number of clusters in the particular iteration. \\
{\em Step 4:} Letting $p_\gamma$ be the dimension of the current model, update $\tau^{-1}$ using \\
$\pi(\tau^{-1}|-) = Ga\bigg(\frac{n+p_\gamma}{2},\frac{1}{2}\left\{ (Y^n-X_{\gamma}\beta_{\gamma})'\Sigma_A^{-1}(Y^n-X_{\gamma}\beta_{\gamma}) +  \frac{1}{g}\beta_{\gamma}'(X_{\gamma}'\Sigma_A^{-1}X_{\gamma})\beta_{\gamma} \right\}\bigg)$.\\
{\em Step 5:} Using the hyper-$g$ prior and the fact that $\frac{g}{1+g}\sim Be(1,a/2-1)$, we can subsequently adopt the griddy Gibbs approach (Ritter and Tanner, 1992) to update $g$.\\
{\em Step 6:} For variable selection, we update $\gamma^j$'s one at a time by computing their posterior inclusion probabilities after marginalizing out $\beta$ and conditional on inclusion indicators for the remaining predictors as well as $g,\tau^{-1}$ and A. Denoting $\gamma(j)$ as the vector of variable inclusion indicators with $\gamma^j=1$, and p$_{\gamma(j)}$ as the vector sum of $\gamma(j)$, we can sample $\gamma^j$ from the Bernoulli conditional posterior distribution with probabilities 
\begin{eqnarray*}
P(\gamma^j=1|-)\propto (1+g)^{-p_{\gamma(j)}/2}\exp\left\{-\frac{\tau^{-1}}{2}\frac{g}{1+g}\bigg(Y^{n'}\Sigma_A^{-1}X_{\gamma(j)}(X_{\gamma(j)}'\Sigma_A^{-1}X_{\gamma(j)})^{-1}X_{\gamma(j)}'\Sigma_A^{-1}Y^n\bigg)\right\}.
\end{eqnarray*}
{\em Step 7:} Set $\left\{\beta_j:\gamma^j=0\right\}=0$ and update $\beta_\gamma=\left\{\beta_j:\gamma^j=1\right\}$ using $\pi(\beta_\gamma|-) = N(\beta_\gamma;E,V)$, \\ where  $V=\bigg(\frac{\tau^{-1}}{g}(X_{\gamma}'\Sigma_A^{-1}X_{\gamma})+\tau^{-1}(X_{\gamma}'X_{\gamma})\bigg)^{-1}$ and $E=V\bigg(\tau^{-1}X'_{\gamma}(Y^n-\alpha)\bigg)$.

\vskip 12pt

{\noindent{\bf 5. \quad SIMULATION STUDY}}\\
We present the results of two simulation studies to demonstrate the utility of Bayesian variable selection in semi-parametric linear models. For the first case (Case I), the truth was generated from a linear model involving ten predictors with coefficients (3 2 -1 0 1.5 1 0 -4 -1.5 0) and a bimodal residual specified by 0.5*N(2.5,1)+0.5*N(-2.5,1). For the second case (Case II), the truth was generated from a normal linear model with the same set of regression coefficients and intercept=1. The covariates were generated independently from uniform(-1,1) distribution. For each case, we generated 20 different replicates for each of the sample sizes 100, 200, 300, 400 and 500, and summarized the results across the replicates.

After generating the data in such a manner, we compared the performance of our method using marginal inclusion probabilities for each predictor (given by P($\beta_j\ne 0|Y^n$), j=$1,\ldots,p$), with the normal linear model having $\beta_\gamma \sim N(0,g\tau^{-1}(\tilde{X}'_{A=1_n,\gamma} \tilde{X}_{A=1_n,\gamma} )^{-1} ) $. This prior on $\beta_{\gamma}$ is a special case of the SLM with A=$1_n$, and is an attempt to assign comparable prior information to both the methods. The replicate averaged marginal inclusion probabilities of each predictor are reported across different sample sizes. We used $Ga(0.1,1)$ prior on the DP precision parameter. Further, we chose a $Be(1,1)$ prior for $\frac{g}{1+g}$ which corresponds to a=4 in the hyper-$g$ prior. For the griddy Gibbs approach, we chose 1,000 equally spaced quantiles from $Be(1,1)$ distribution. We made 50,000 runs with a burn in of 5,000. 

As the sample size increases, it is interesting to see how the marginal probabilities of inclusion for important predictors and the marginal probabilities of exclusion for unimportant predictors change. The marginal inclusion probabilities under both the methods were 1.00 for $\beta_1=3$ and $\beta_8= -4$ for all the sample sizes. For the remaining predictors, the plots of the marginal inclusion probabilities over different sample sizes are presented (Fig 1 and Fig 2), as a comparison between the two methods. These plots depict a faster rate of increase of the marginal inclusion probabilities of the important predictors for the semi-parametric Bayes method when the true residuals are non-Gaussian and a similar rate of increase under both methods when the true residuals are Gaussian. In contrast, for the unimportant predictors the exclusion probabilities converge to one slowly for both the methods, reflecting the well known tendency for slower accumulation of evidence in favor of the true null.

To get a closer look when the true residual is non-Gaussian (Case I), we present the results for the sample size 100. As a comparison, we also present regression estimates under the lasso (Tibshirani, 1996) and elastic net (Zou and Hastie, 2005), using the GLMNET package in R with default settings. The average mean square error for out of sample prediction for a test sample size of 25 under the semi-parametric linear model was 7.7 compared to 15.4 under the normal linear model, implying a 50\% reduction. The average out of sample MSE were 7.68 for lasso (L1) and 7.65 for elastic net (EL). Out of the 20 different replicates generated with sample size 100, the normal linear model (NLM) chose the wrong subset of predictors 16 times under the median probability model, whereas the semi-parametric linear model (SLM) made incorrect variable selection decisions for 3 out of 20 replicates. The computation time for SLM per iteration was marginally slower than NLM, with the difference inreasing as the number of clusters increase. The mixing for the fixed effects was good under both the methods. The results for SLM do not appear to be sensitive to the hyper-parameters in $\pi(m)$, but are mildly sensitive to hyper-parameters in $\pi(g)$ for n=100.

Table 1 summarizes results for the model averaged regression estimates ($\hat{\beta}$), including 95\% pointwise credible intervals (C.I.) and marginal inclusion probabilities (MIP). The SLM and NLM correctly identify the important as well as unimportant predictors. In general, the L1 and EL results seem to be unstable with the coefficients shrunk to 0 varying over different replicates. As a result, the replicate averaged estimates for L1 and EL in Table 1 lead to inaccurate estimates of $\beta_4$ and $\beta_7$. For the estimation of the fixed effects, the MSE around the true $\beta$ ($\frac{||\hat{\beta}-\beta_{True} ||_2}{p}$) was 0.015 for SLM, 0.084 for NLM, 0.047 for elastic net and 0.047 for lasso. Thus, the SLM results in more accurate estimates with narrower credible intervals. From the results, it is clear that when the true residual is non-Gaussian, the SLM not only has a more desirable performance in variable selection and estimation of regression coefficients, it also has a superior out of sample predictive performance as compared to NLM.  
\vskip 12pt

{\noindent{\bf 6. \quad APPLICATION TO DIABETES DATA}}\\
The prevalence of diabetes in the United States is expected to more than double to 48 million people by 2050 (Mokdad et. al., 2001). Previous medical studies have suggested that Diabetes Mellitus type II (DM II) or adult onset diabetes could be associated with high levels of total cholesterol (Brunham et. al., 2007) and obesity (often characterized by BMI and waist to hip ratio) (Schmidt et. al., 1992), as well as hypertension (indicated by a high systolic or diastolic blood pressure or both) which is twice as prevalent in diabetics compared to non-diabetic individuals (Epstein and Sowers, 1992). However, most of these results rely on informal treatment of data and lack rigorous statistical analysis to support their conclusions. 

We develop a comprehensive variable selection strategy for indicators of DM II based on data obtained from Department of Biostatistics, Vanderbilt University website, involving a diabetes study for African-Americans. Our primary focus is to discover important indicators of DM II by modeling the continuous outcome, glycosylated hemoglobin ($>7 mg/dL$  indicates a positive diagnosis of diabetes) based on predictors such as total cholesterol (TC), stabilized glucose (SG), high density lipoprotein (HDL), age, gender, body mass index (BMI) indicator (overweight and obese with normal as baseline), systolic blood pressure (SBP), diastolic blood pressure (DBP), waist to hip ratio (WHR) and postprandial time indicator (PPT) (0/1 depending on whether the blood was drawn within 2 hours of a meal). In addition to the factors already noted above (total cholesterol, obesity and hypertension) for DM II, we note that lower levels of HDL have been known to be associated with insulin resistance syndrome (often considered a precursor of DM II with a conversion rate around 30\%), and further we also expect PPT to be a significant indicator as blood sugar levels are high up to 2 hours after a meal.

After trimming the records containing missing values, the data consisted of 365 subjects which was split into multiple training and test samples of sizes 330 and 35 respectively. The replicate averaged fixed effects estimates (multiplied by 100) for the SLM, NLM, lasso (L1) and elastic net (EL) are presented in Table 2, along with the marginal inclusion probabilities (MIP) for the SLM and the NLM. We also evaluate the out of sample predictive performance for each training-test split using predictive MSE for SLM, NLM, L1 and EL in Table 3, and additionally provide coverage (COV) and width (CIW) of 95\% pointwise credible intervals for SLM and NLM. The same values of hyper-parameters were used as in section 5 for SLM and NLM. For each replicate, we randomized the initial starting points and made 100,000 runs for SLM (burn in = 20,000) and 50,000 runs for NLM (burn in = 5,000). 

It is interesting to note from Table 2 that SLM tells a quite different story compared to the NLM in terms of variable selection. In particular, while both the models successfully identify total cholesterol, stabilized glucose and postprandial time as important predictors, it is only the SLM which identifies systolic hypertension (MIP = 0.77) and waist to hip ratio (MIP = 0.98) as important positively associated indicators, whereas NLM fails to identify these factors (MIP = 0.18 for SBP and 0.17 for WHR) and instead throws in age (MIP = 0.77) as an important predictor. Further, SLM points to a more significant  negative association with HDL (MIP=0.64) as compared to NLM (MIP=0.52). For both the methods, the marginal inclusion probabilities for BMI (overweight and obese) were low, which could potentially be attributed to adjusting for the other factors such as waist to hip ratio. The lasso and elastic net include all predictors except DBP, and hence produce an overly complex model. 

Variable selection in this application is clearly influenced by the assumptions on the residual density, with the nonparametric residual density providing a more realistic characterization that should lead to a more accurate selection of the important predictors.  Figure 3 show an estimate of the residual density obtained from the SLM analysis, suggesting a unimodal right skewed density with a heavy right tail.  The SLM results suggest that a mixture of two Gaussians provides an adequate characterization of this density. The computation time for SLM is only marginally slower than NLM, and in addition SLM exhibits good mixing for most of the fixed effects (Table 4). These results are robust to SSVS starting points, and consistency in the results across training-test splits also indirectly suggests adequate computational efficiency of SSVS. 

In terms of out of sample predictive MSE (Table 3), none of the models is a clear winner, with the relative performance varying across training-test splits. The MSE's for lasso and elastic net are very similar to NLM, except for the second test sample where they have the lowest MSE. Overall, the NLM has narrower 95\% pointwise credible intervals compared to SLM, often resulting in poorer coverage. Thus, in conclusion, although the competitors yield comparable out of sample predictive performance, it is only the SLM which succeeds in choosing the most reasonable model for DM II, consistent with previous medical evidence. 
		
\vskip 12pt

{\noindent{\bf 7. \quad DISCUSSION}}\\
We develop mixtures of semi-parametric $g$-priors for linear models with non-parametric residuals characterized by DP mixtures of Gaussians. The proposed method addresses the often encountered issue of non-Gaussianity of residuals in variable selection settings, and has attractive asymptotic justifications such as Bayes factor and variable selection consistency involving fixed $p$ as well as $p>n$ (under some restrictions on the model space). Further, the method is essentially no more difficult to implement than SSVS for normal linear models and can lead to substantially different conclusions, as illustrated in the diabetes application.  The general topic of semi- and nonparametric Bayesian model selection is understudied and we hope that this work stimulates additional research of this type in broader model classes, such as for generalized linear models and nonparametric regression. 
\vskip 12pt

{\noindent{\bf 7. \quad ACKNOWLEDGEMENTS}}\\
This work was support by Award Number R01ES017240 from the National Institute of Environmental Health Sciences.  The content is solely the responsibility of the authors and does not necessarily represent the official views of the National Institute of Environmental Health Sciences or the National Institutes of Health.
\vskip 12 pt 

{\noindent {\quad APPENDIX: PROOF OF RESULTS}}\\
\textbf{Proof of Theorem I}: Using similar methods as in the proof of theorem 2 in Guo and Speckman (2009), it can be shown that conditional on A and assumptions \emph{(A3)} and \emph{(A4)}, the upper and lower bounds of 
$I_1=\int_0^\infty (1+g)^{-p_1/2}\left[1-\frac{g}{1+g}\tilde{R}_{A,1}^2\right]^{-n/2}\pi(dg)$ are 
\begin{eqnarray*}
I_1 &\le& \bigg(\frac{p_1+2k_u}{n-p_1-2k_u}\bigg)^{p_1/2+k_u}\bigg(\frac{1-\tilde{R}_{A,1}^2}{\tilde{R}_{A,1}^2}\bigg)^{p_1/2+k_u}\bigg(\frac{n}{n-p_1-2k_u}\bigg)^{-n/2}\bigg(1-\tilde{R}_{A,1}^2\bigg)^{-n/2} \\
&\approx& \bigg(\frac{p_1+2k_u}{n-p_1-2k_u}\bigg)^{p_1/2+k_u}\bigg(\frac{1-\tilde{R}_{A,1}^2}{\tilde{R}_{A,1}^2}\bigg)^{p_1/2+k_u}\bigg(1-\tilde{R}_{A,1}^2\bigg)^{-n/2} = U_{A,1}(n), 
\end{eqnarray*}
and $I_1 \ge n^{-p_1/2-k}\bigg(1-\tilde{R}_{A,1}^2\bigg)^{-n/2} = L_{A,1}(n)$. Similarly,  
\begin{eqnarray*}
 L_{A,2}(n) &\le& I_2 = \int_0^\infty (1+g)^{-p_2/2}\left[1-\frac{g}{1+g}\tilde{R}_{A,2}^2\right]^{-n/2}\pi(dg) \le U_{A,2}(n). 
\end{eqnarray*}
Therefore, BF$^n_{21,A} \le \frac{U_{A,2}(n)}{L_{A,1}(n)} $
\begin{eqnarray}  = \bigg(\frac{p_2+2k_u}{n-p_2-2k_u}\bigg)^{p_2/2+k_u}\bigg(\frac{1-\tilde{R}_{A,2}^2}{\tilde{R}_{A,2}^2}\bigg)^{p_2/2+k_u}\bigg(1-\tilde{R}_{A,2}^2\bigg)^{-n/2}/\bigg( n^{-p_1/2-k} (1-\tilde{R}_{A,1}^2)^{-n/2} \bigg) . \label{eq:bddBF}
\end{eqnarray}
Case (I): For fixed $p$, directly from the proof of Theorem 3 in Guo and Speckman (2009) 
\begin{eqnarray}
\mbox{BF}^n_{21,A}&\le& \zeta(A,n) = n^{\frac{p_1-p_2}{2}+k-k_u}\bigg(\frac{1-\tilde{R}_{A,2}^2}{1-\tilde{R}_{A,1}^2}\bigg)^{-n/2} \to 0 \mbox{ under } \mathcal{M}_1 \mbox{ for all A}\in \mathcal{C}_\infty.  \label{eq:zeta} 
\end{eqnarray}
\begin{eqnarray}
\mbox {Further, } \quad \mbox{BF}^n_{21,A}&\le& \zeta(A,n)\Rightarrow L(Y^n|A,\mathcal{M}_2) \le \zeta(A,n)L(Y^n|A,\mathcal{M}_1) \nonumber \\
\Rightarrow
L(Y^n|\mathcal{M}_2) &\le& \sum_{A_l\in \mathcal{C}_n}w_l\zeta(A_l,n)L(Y^n|A_l,\mathcal{M}_1) \le \mbox{max}_{A \in \mathcal{C}_n} \zeta(A,n) L(Y^n|\mathcal{M}_1). \label{eq:bddzeta} 
\end{eqnarray}
In the limiting sense as $n\to \infty$, the maximum in the upper bound in (\ref{eq:bddzeta}) is computed over A$\in \mathcal{C}_\infty$. From (\ref{eq:zeta}), $\zeta(A,n)\to 0$ under $\mathcal{M}_1$ for all A$\in \mathcal{C}_\infty$ implies max$_{A \in \mathcal{C}_\infty} \zeta(A,n)\to 0$. Dividing both sides of (\ref{eq:bddzeta}) by L($Y^n|\mathcal{M}_1$), this implies BF$^n_{21}\to$0 under $\mathcal{M}_1$. Further, the mode of convergence of BF$^n_{21}$ is the same as BF$^n_{A,21}$, and the rest follows from the proof of Theorem 3 in Guo and Speckman (2009). 

Case (II): For increasing model dimensions $p_1=O(n^{a_1})$ and $p_2=O(n^{a_2})$ with $0\le a_1<a_2<1$, for $g \sim \pi(g)$ we will only assume \emph{(A3)} so that k$_u=0$. We have using (\ref{eq:bddBF})
\begin{eqnarray}
BF^n_{21,A} \le n^{p_1/2-(1-a_2)p_2/2+k}\bigg(\frac{1-\tilde{R}_{A,2}^2}{\tilde{R}_{A,2}^2}\bigg)^{p_2/2}\bigg(\frac{1-\tilde{R}_{A,2}^2}{1-\tilde{R}_{A,1}^2}\bigg)^{-n/2}.  \label {eq:boundBF}
\end{eqnarray}
Let us consider the following cases under $0 \le a_1<a_2<1$. \\
Case C1: $\mathcal{M}_1 \subset \mathcal{M}_2$. We have  $Q_j=\tau(Z_A'Z_A - Z_A'\tilde{H}_{A,j}Z_A) \sim \chi^2_{n-p_j}(0)$, j=1,2, and $Q_1-Q_2 = \tau \bigg(Z_A'(\tilde{H}_{A,2} - \tilde{H}_{A,1})Z_A\bigg) \sim \chi^2_{p_1-p_2}(0)$. Using Lemma 1 of Guo et. al. (2009), 
\begin{eqnarray*}
\frac{1-\tilde{R}_{A,2}^2}{1-\tilde{R}_{A,1}^2}= \frac{Z_A'Z_A - Z_A'\tilde{H}_{A,2}Z_A}{Z_A'Z_A - Z_A'\tilde{H}_{A,1}Z_A} = \frac{Q_2}{Q_1}= 1 - \frac{(Q_1-Q_2)/(p_2-p_1)}{Q_1/(n-p_1)}\frac{p_2-p_1}{n-p_1}\stackrel{a.s.}{\rightarrow} 1. 
\end{eqnarray*}
Moreover $\bigg(\frac{1-\tilde{R}_{A,2}^2}{\tilde{R}_{A,2}^2}\bigg)\stackrel{a.s.}{\rightarrow} \bigg(\frac{\tau^{-1}}{b_{A,1}}\bigg)$ under $\mathcal{M}_1$, which implies that $\bigg(\frac{1-\tilde{R}_{A,2}^2}{\tilde{R}_{A,2}^2}\bigg)^{p_2/2}$ blows up at a rate strictly slower than the rate at which $n^{p_1/2-(1-a_2)p_2/2+k} \to 0$. This implies that BF$^n_{21,A}\stackrel{a.s.}{\rightarrow}0$ under $\mathcal{M}_1$, for all A$\in \mathcal{C}_\infty$.  \\
Case C2: $\mathcal{M}_1\not\subseteq\mathcal{M}_2$. Using Lemma 1, 
\begin{eqnarray*}
\frac{1-\tilde{R}_{A,2}^2}{\tilde{R}_{A,2}^2} \stackrel{a.s.}{\rightarrow} \frac{\tau^{-1}+b_{A,1}-b_{A,2}}{b_{A,2}}, \quad \frac{1-\tilde{R}_{A,1}^2}{1-\tilde{R}_{A,2}^2} \stackrel{a.s.}{\rightarrow} \frac{\tau^{-1}}{\tau^{-1}+b_{A,1}-b_{A,2}}, \mbox{ under }\mathcal{M}_1.
\end{eqnarray*} 
For fixed $\tau^{-1}$ and $b_{A,2}>0$ (under (\emph{$\tilde{A}2$})),  $\bigg(\frac{1-\tilde{R}_{A,2}^2}{\tilde{R}_{A,2}^2}\bigg)^{p_2/2}\bigg(\frac{1-\tilde{R}_{A,2}^2}{1-\tilde{R}_{A,1}^2}\bigg)^{-n/2}\stackrel{a.s.}{\rightarrow}0$. In addition, we have $p_1-(1-a_2)p_2+k<0$ for $0\le a_1<a_2<1$, which implies  BF$^n_{21,A}\stackrel{a.s.}{\rightarrow}0$ under $\mathcal{M}_1$. \\
Subsequently using similar arguments as in Case (I), BF$^n_{21} \stackrel{a.s.}{\rightarrow}0$ under $\mathcal{M}_1$ for both C1, C2. 

\textbf{Proof of Theorem II}: For fixed $p$, the result is trivial to prove using Theorem I. For increasing $p$ ($\ge p_1$),  $\pi^n(\gamma_{j_{l}}) \propto  2^{-p_j/2}I[p_j\le O(n^{a_1})] +  N_j^{-1}I[O(n^{a_1})<p_j\le (n-1)\wedge p]$, $l=1,\ldots,N_j$, $N_j = \binom{p}{p_j}$. Let BF$_{j_{l}1}$= Bayes factor between models $\gamma_{j_l}$ and $\mathcal{M}_1$ and H$^n(a_1):=\left\{j:O(n^{a_1})<p_j\le(n-1)\wedge p\right\}$. Then  
$ P(\mathcal{M}_1|Y^n) = [ 1 + 2^{p_1/2}\sum_{j:p_j\in H^n(a_1)} \sum_{l=1}^{N_j} N_j^{-1}BF^n_{j_{l}1}] ^{-1}.$
From the preceeding proof of Theorem I, the upper bound for $\left\{\mbox{BF}^n_{j_{l}1}:O(n^{a_1})<p_j<(n-1)\wedge p\right\}$ for large $n$ is $\bar{U}^n_j$, given by (a) $\bar{U}^n_j \approx$ $\kappa^{p_j/2}n^{-(1-a_j)p_j/2 + p_1/2+k}$ for  $0<\kappa<\infty$, for the nested
case (b) $\bar{U}^n_j\le n^{-(1-a_j)p_j/2 + p_1/2+k}$, for non-nested case, with k$\ge 0$. Therefore
\begin{eqnarray*}
P(\mathcal{M}_1|Y^n) &\ge& [1+2^{p_1/2}\sum_{j:p_j\in H^n(a_1)} \sum_{l=1}^{N_j} \bar{U}^n_j/N_j]^{-1} = [1+2^{p_1/2}\sum_{j:p_j\in H^n(a_1)} \bar{U}^n_j]^{-1} \\
&\ge& [1+ n 2^{p_1/2} \mbox{max}_{j:p_j\in H^n(a_1)} \bar{U}^n_j]^{-1} \to 1 \mbox{ as } n\to \infty, \mbox{ using (a), (b), and Theorem I}.
\end{eqnarray*}
Further, the mode of convergence under $\mathcal{M}_1$ is the same as that of the associated conditional Bayes factors.

\vskip 20 pt

\begin{small}
\textbf{Table 1: Results for Case I when n=100.} SLM: Semi-parametric linear model, NLM: Normal linear model, L1: Lasso, EL: Elastic Net. MIP: Marginal Inclusion Probability.
\begin{center}
\begin{tabular}{|c c c c c c c|}
\hline 
$\beta_{True}$ &MIP$_{SLM}$ &MIP$_{NLM}$ &$\hat{\beta}_{SLM}$ &$\hat{\beta}_{NLM}$&$\hat{\beta}_{L1}$ &$\hat{\beta}_{EL}$\\
\hline
3 &0.99 &0.99 &2.89 (2.22,3.44) &2.74 (1.73,3.57) &2.93  &2.92 \\
2 &0.98 &0.94 &1.84 (1.37,2.55) &1.65 (1.21,3.28) &1.61  &1.61     \\
-1 &0.89 &0.59 &-0.84 (-1.47,-0.32) &-0.59 (-1.95,-0.03) &-1.18  &-1.17    \\
0 &0.19 &0.25 &-0.02 (-0.41,0.33) &-0.01 (-0.52,0.45) & 0.21 & 0.21 \\
1.5 &0.97 &0.85 &1.39 (0.89,2.08) &1.24 (1.25,3.12) &1.9 &1.89 \\
1 &0.91 &0.58 &0.83 (0.27,1.42) &0.61 (-0.12,1.11) &0.92 &0.92 \\
0 &0.24 &0.28 &-0.05 (-0.45,0.28) &-0.06 (-0.72,0.29) &-0.22 &-0.23 \\
-4 &1.00 &1.00 &-3.86 (-4.39,-3.16)  &-3.64 (-4.42,-2.52) &-4.05 & -4.04  \\
-1.5 &0.95 &0.77 &-1.34 (-2.03,-0.81)  &-1.14 (-2.11,-0.01) &-1.33  & -1.34\\
0 &0.21 &0.30 &-0.02 (-0.41,0.31) &-0.04 (-0.59,0.41) &0.07  & 0.07  \\
\hline
\end{tabular}
\end{center}
\end{small}  
 \newpage
 
\begin{small}
\textbf{Table 2: Fixed effects (times 100) and marginal inclusion probabilities (MIP).}\\
SLM: Semi-parametric linear model, NLM: Normal linear model, L1: Lasso, EL: Elastic Net.
\begin{center}
\begin{tabular}{|l l l c c c c|}
\hline 
Predictor &$\hat{\beta}_{SLM}$ &$\hat{\beta}_{NLM}$  &$\hat{\beta}_{L1}$ &$\hat{\beta}_{EL}$ &MIP$_{SLM}$ &MIP$_{NLM}$  \\
\hline
TC   & 0.48(0.12,    0.84)           & 0.63( 0.17,    1.07)            &0.75             & 0.75     &0.98   &0.99 \\
   
SG   & 2.16 (1.81,    2.52)          & 2.84 ( 2.52,    3.17)           &2.83             & 2.82     &1.00   &1.00 \\
     
HDL  &-0.48 (-1.33,    0.02)          &-0.52 (-1.66,    0.03)          &-1.02            &-1.02     &0.64   &0.52  \\
     
Age  & 0.51 ( 0,       1.64)          & 1.30 ( 0.01,    2.56)          &1.19             & 1.19     &0.35   &0.77  \\
     
Gender   &-3.05 (-28.83, 6.98)        & -1.90 (-26.34,    7.14)        &-19.66           &-19.81    &0.21   &0.16  \\
        
BMI(overwt)   & 0.82 (-7.54, 18.52)         & 2.04 (-5.30, 29.59)         &4.33   & 4.27     &0.15   &0.17  \\
                  
BMI(obese)   & 0 (-13.12,   13.00)                & -1.37 (-24.83, 9.13)      &-14.88 &-15.03     &0.15   &0.15  \\
            
SBP   & 0.45 (-0.02,    1.24)            & 0.04 (-0.19,    0.71)              &0.25   &0.25     &0.77  &0.18  \\
      
DBP   &-0.03 (-0.94,    0.61)             & 0 (-0.58,    0.56)                &0.018  &0.017     &0.19  &0.14  \\
     
WHR   &211.74 (40.02,    361.41)           & 4.72(-53.12,  102.12)            &90.47  &91.53     &0.98   &0.17  \\
     
PPT   &20.62(0,    56.13)               & 32.13 (0,   75.89)                  &47.31  &47.32     &0.62   &0.77  \\
      
\hline
\end{tabular}
\end{center}
\end{small}

\begin{small}
\textbf{Table 3: Out of Sample Prediction.} (Cov: 95\% C.I. coverage, CIW: 95\% C.I. width) 
\begin{center}
\begin{tabular}{|l c c c c c c c c|}
\hline 
&MSE$_{SLM}$ &Cov(SLM) &CIW(SLM) &MSE$_{SLM}$ &Cov(NLM) &CIW(NLM) &MSE$_{L1}$ &MSE$_{EL}$\\
\hline
Sample 1   & 1.27  & 97.14   & 6.94   & 1.23  & 97.14   & 5.91  & 1.36 & 1.24 \\
Sample 2   & 4.67  & 94.28   & 6.23   & 4.67  & 91.42   & 5.40  & 1.21 & 1.20 \\
Sample 3   & 1.55  &100.00   & 6.83   & 1.78  & 94.28   & 5.82  & 1.75 &1.75 \\
Sample 4   & 1.22  & 97.14   & 6.77   & 1.26  & 97.14   & 5.91  & 1.24 &1.23 \\
Sample 5   & 1.42  &100.00   & 6.79   & 1.16  &100.00   & 5.92  & 1.17 &1.18 \\
Sample 6   & 1.46  &100.00   & 6.79   & 1.43  & 97.14   & 5.90  & 1.52 &1.52 \\
Sample 7   & 3.70  & 91.42   & 6.47   & 3.40  & 91.42   & 5.59  & 3.38 &3.38\\
Sample 8   & 1.24  &100.00   & 6.79   & 1.50  & 97.14   & 5.87  & 1.54 &1.53 \\
\hline
\end{tabular}
\end{center}
\end{small}

\begin{small}
\textbf{Table 4: Auto-correlations across lags for fixed effects.} 
\begin{center}
\begin{tabular}{|l| l l l l l|}
\hline 
Predictor     & Lag 1 & Lag 5 & Lag 10 & Lag 25 & Lag 50 \\
\hline
     &SLM \quad NLM      & SLM \quad NLM      & SLM \quad NLM      & SLM \quad NLM      & SLM \quad NLM \\
TC   &0.22\quad 0.18 & 0.113\quad 0.194 & 0.073 \quad 0.159 & 0.032\quad 0.111 & 0.013\quad 0.059 \\
SG   &0.59\quad 0.06 & 0.386\quad 0.038 & 0.285\quad 0.022 & 0.14\quad 0.009 & 0.06\quad 0.016 \\
HDL   &0.19\quad 0.02 & 0.081\quad 0.012 & 0.041\quad 0.013 & 0.01\quad 0.021 & 0.0005\quad-0.006 \\
Age   &0.21\quad 0.04 & 0.072\quad 0.009 & 0.053\quad-0.0001 & 0.025\quad 0.006 & 0.007\quad-0.014 \\
Gender   &0.06\quad-0.007 & 0.030\quad 0.0003 & 0.013\quad-0.006 & 0.009\quad-0.014 & 0.005\quad 0.019 \\
BMI(overwt)   &0.02\quad-0.002 & 0.01\quad -0.006 & 0.006\quad 0.013 &-0.006\quad 0.009 & 0.0014\quad 0.018 \\
BMI(obese)   &0.02\quad 0.002 & 0.017\quad 0.004 & 0.004\quad 0.018 & 0.007\quad-0.003 & 0.000 \quad 0.000 \\
SBP   &0.29\quad 0.0711 & 0.137\quad 0.019 & 0.096\quad 0.007 & 0.047\quad 0.03 & 0.014\quad 0.022 \\
DBP   &0.07\quad 0.0239 & 0.021\quad 0.019 & 0.019\quad 0.031 & 0.009\quad-0.003 & 0.004\quad-0.012 \\
WHR  &0.44\quad 0.0642 & 0.353\quad 0.043 & 0.321\quad 0.061 & 0.251\quad 0.06 & 0.186\quad-0.003 \\
PPT  &0.22\quad 0.0600 & 0.118\quad 0.047 & 0.068\quad 0.045 & 0.015\quad 0.004 &-0.002\quad 0.019 \\
\hline
\end{tabular}
\end{center}
\end{small}
 
\begin{figure}
		\mbox{\includegraphics[height=2.0in, width=0.45\textwidth]{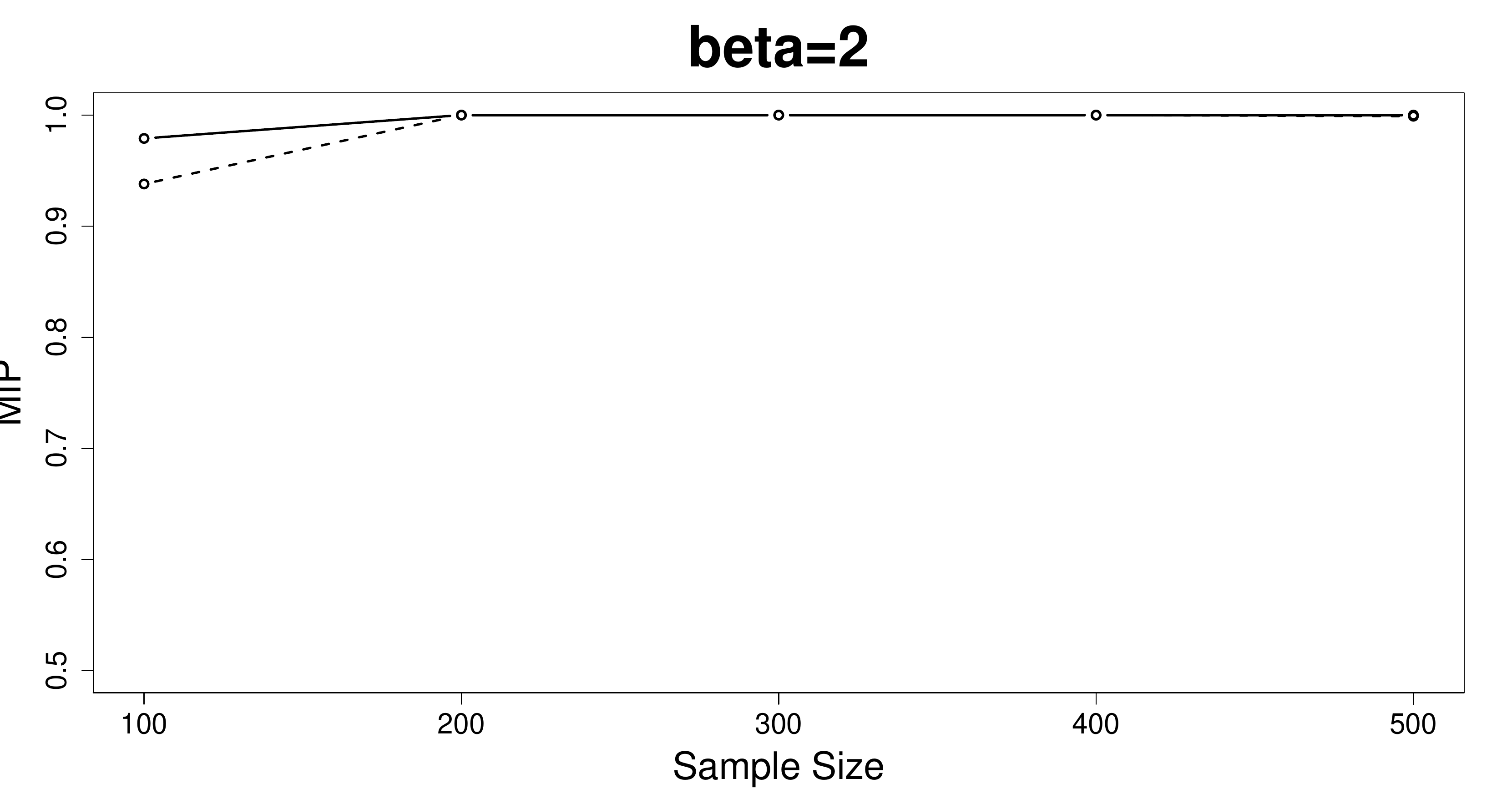}} \quad \mbox{\includegraphics[height=2.0in, width=0.45\textwidth]{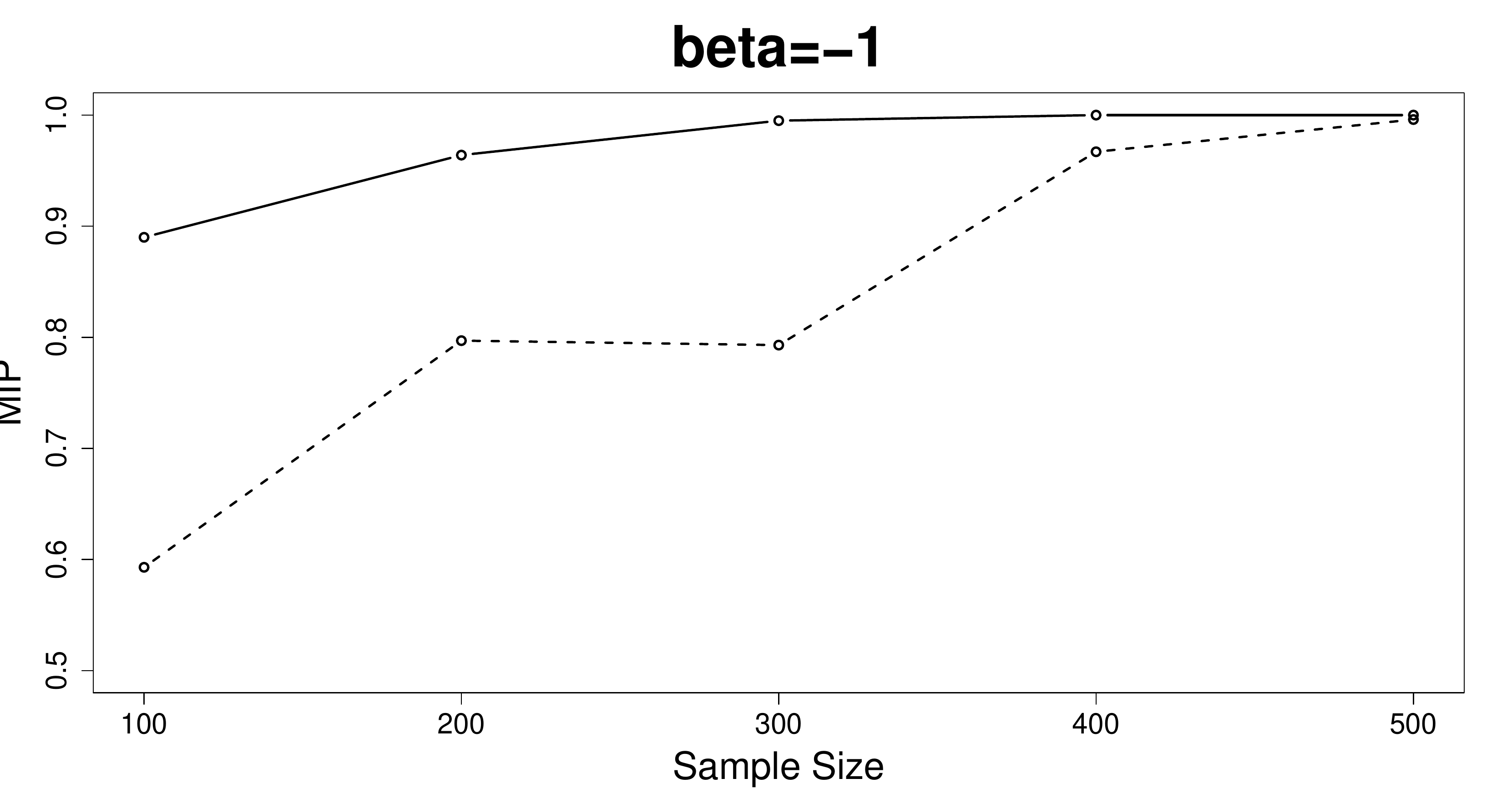}}\\
		\mbox{\includegraphics[height=2.0in, width=0.45\textwidth]{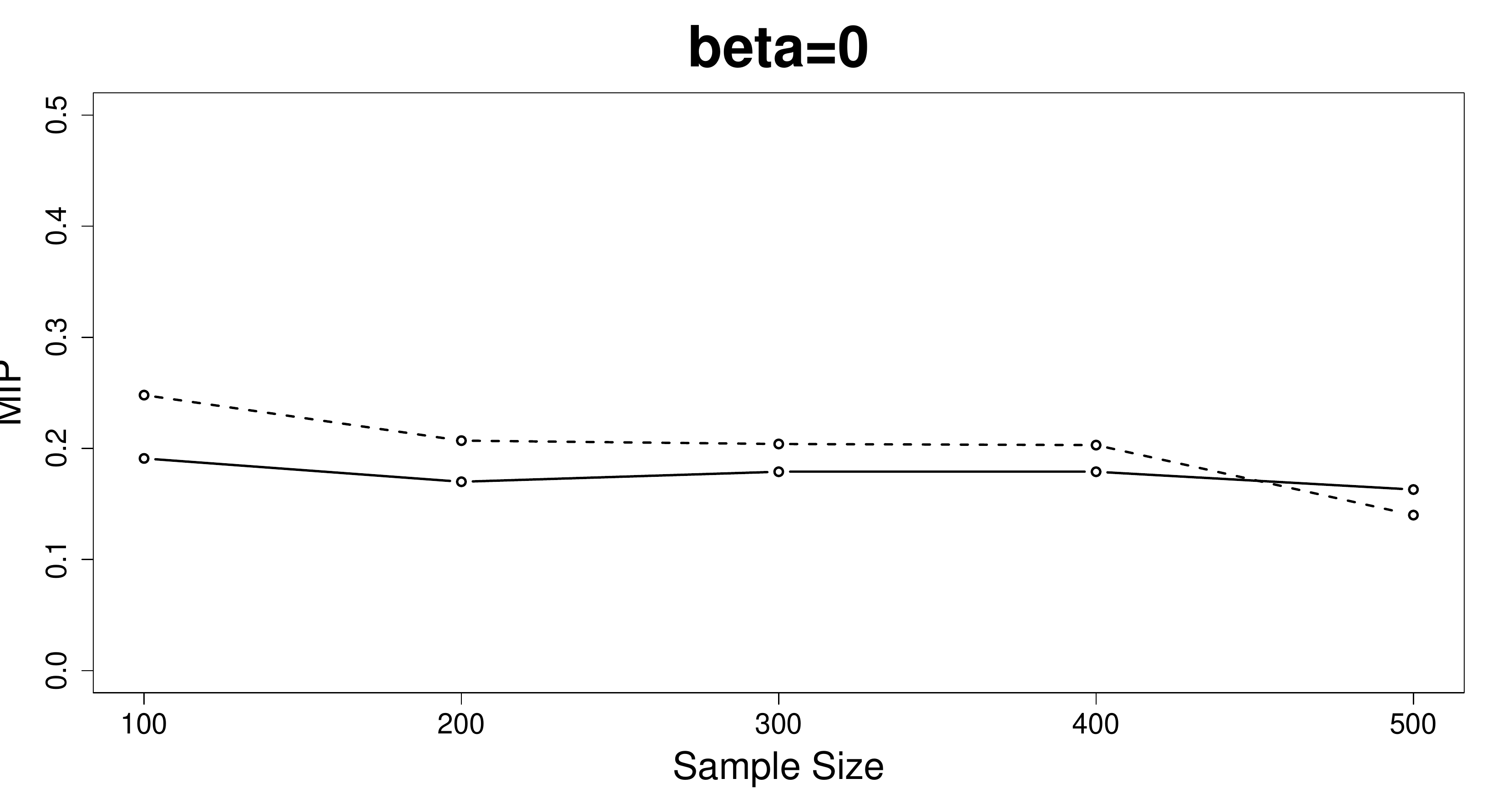}} \quad 	\mbox{\includegraphics[height=2.0in, width=0.45\textwidth]{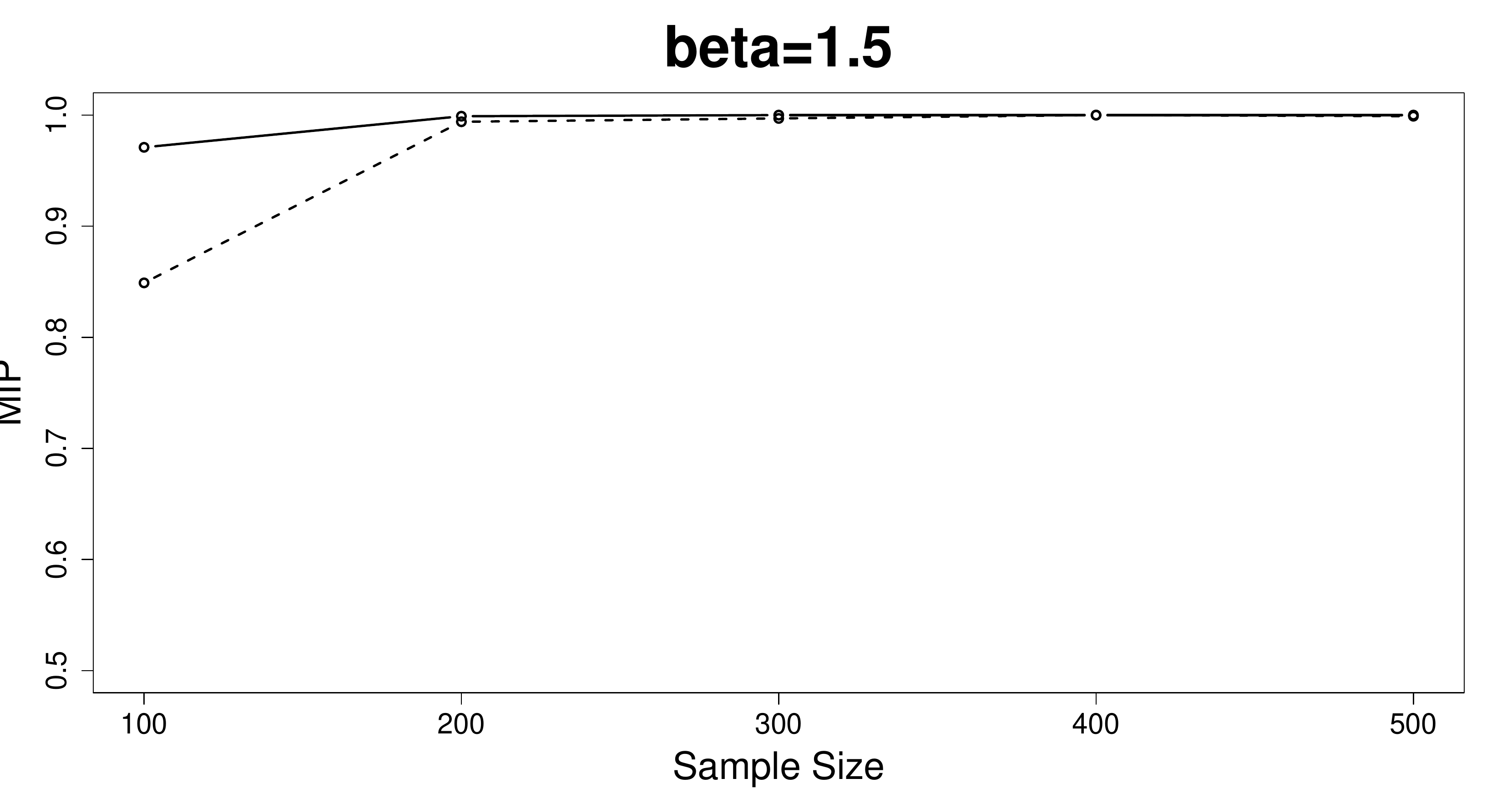}}\\
			\mbox{\includegraphics[height=2.0in, width=0.45\textwidth]{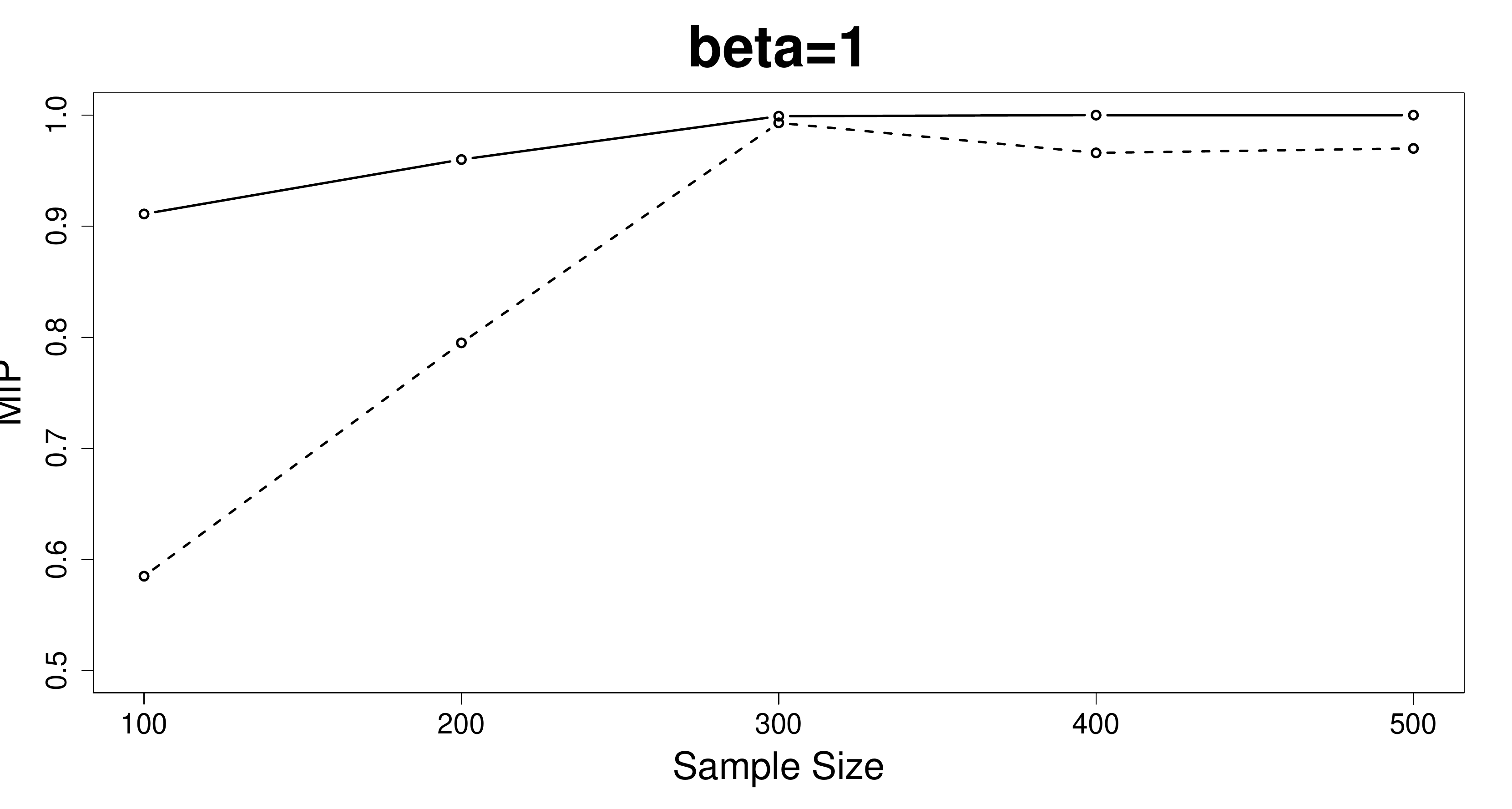}} \quad \mbox{\includegraphics[height=2.0in, width=0.45\textwidth]{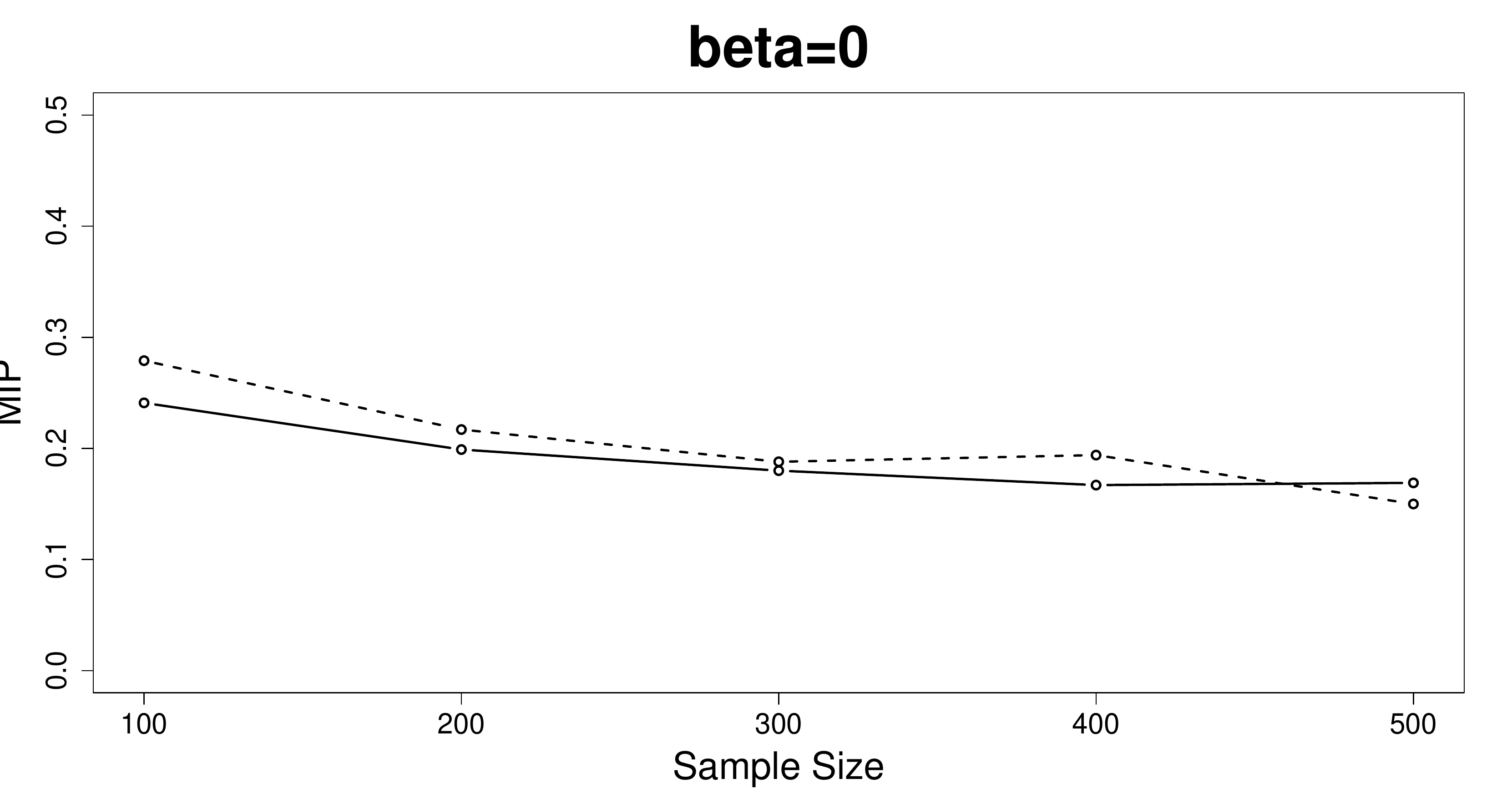}}\\
		\mbox{\includegraphics[height=2.0in, width=0.45\textwidth]{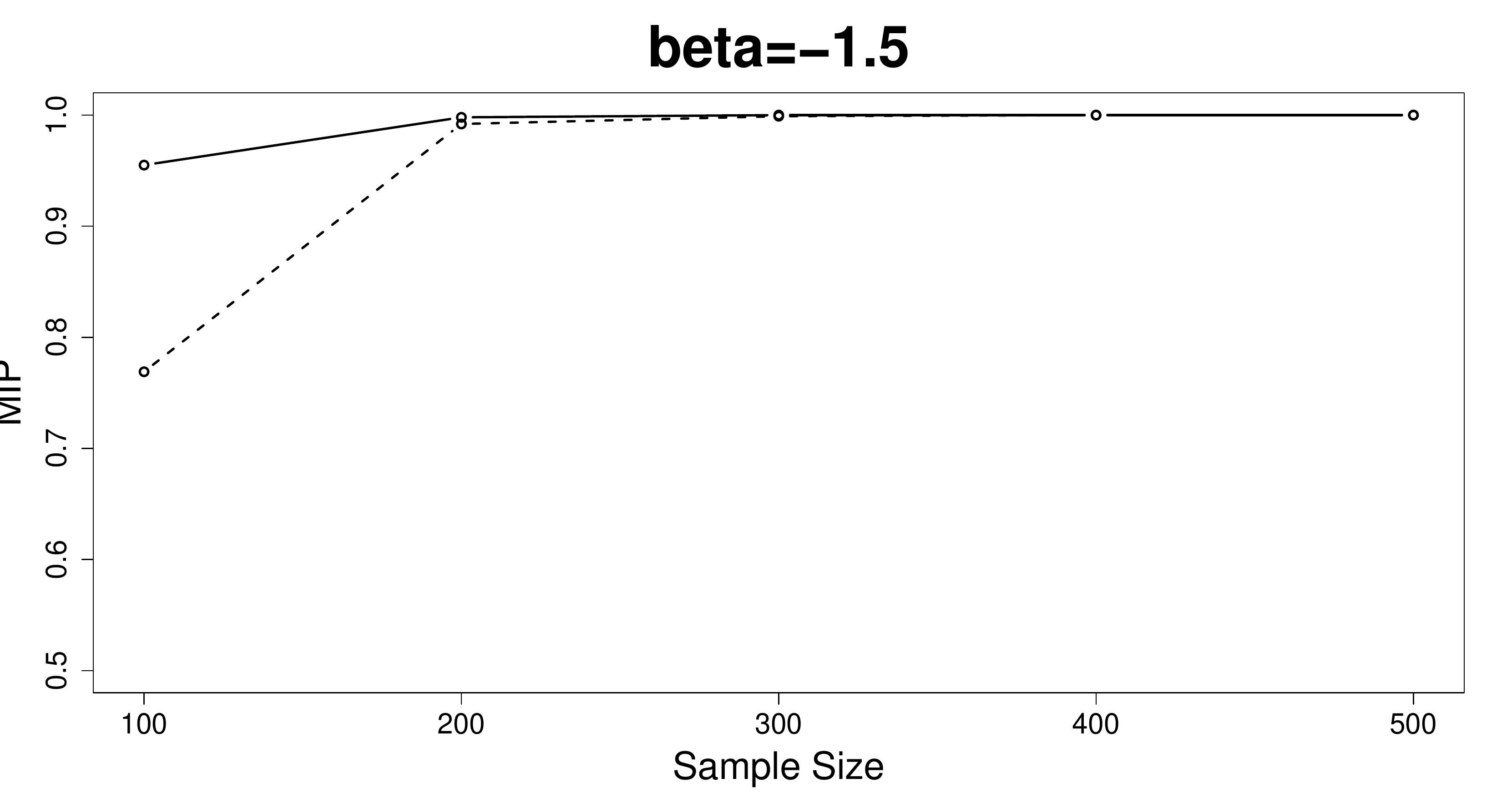}} \quad 	\mbox{\includegraphics[height=2.0in, width=0.45\textwidth]{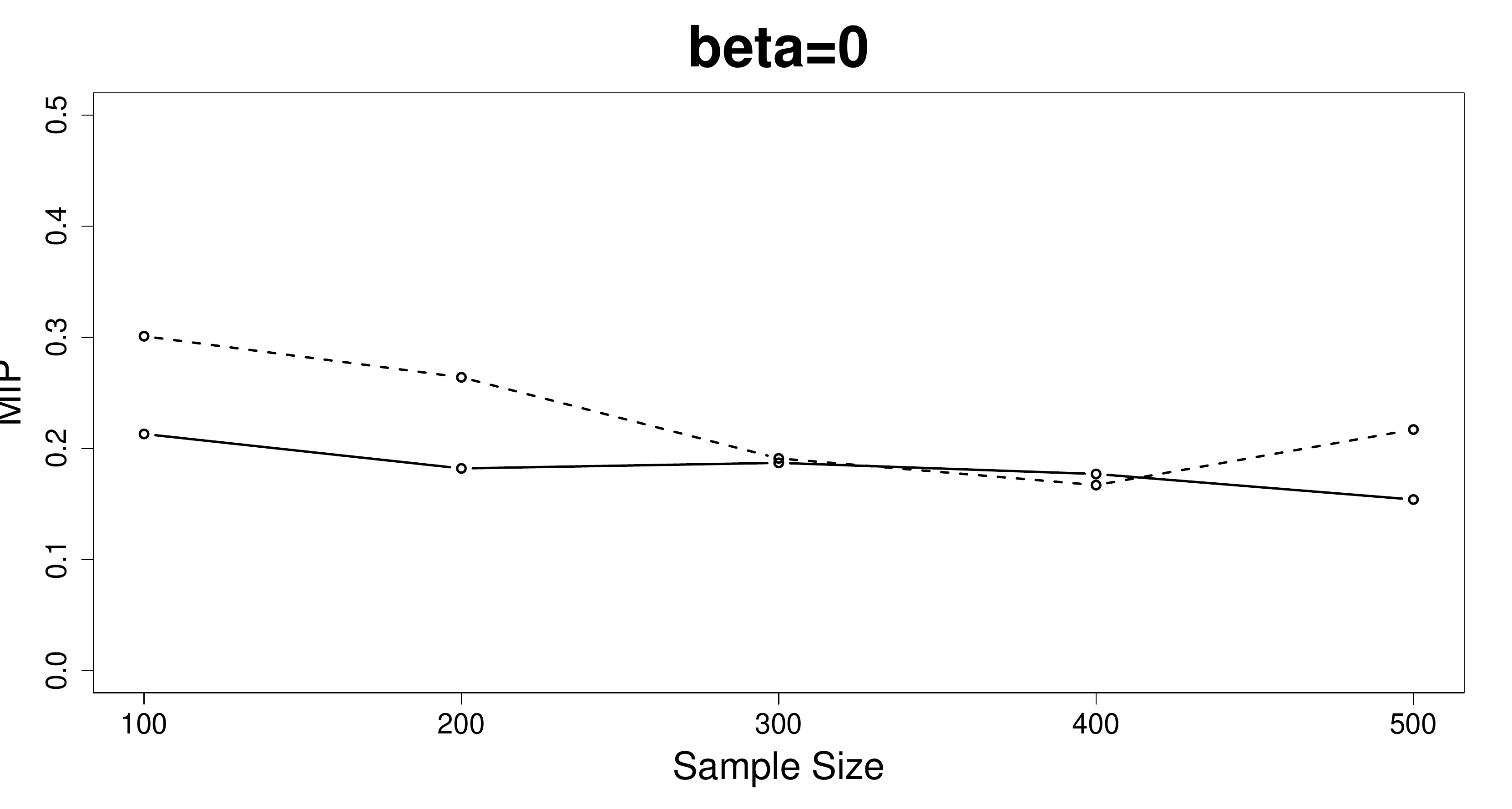}}
			\caption{\small{Marginal Inclusion Probabilities (MIP): Truth generated from bimodal residual. Solid lines - Semi-parametric Linear Model, dashed lines - Non-parametric Linear Model.}}
\end{figure}     
\begin{figure}
		\mbox{\includegraphics[height=2.0in, width=0.45\textwidth]{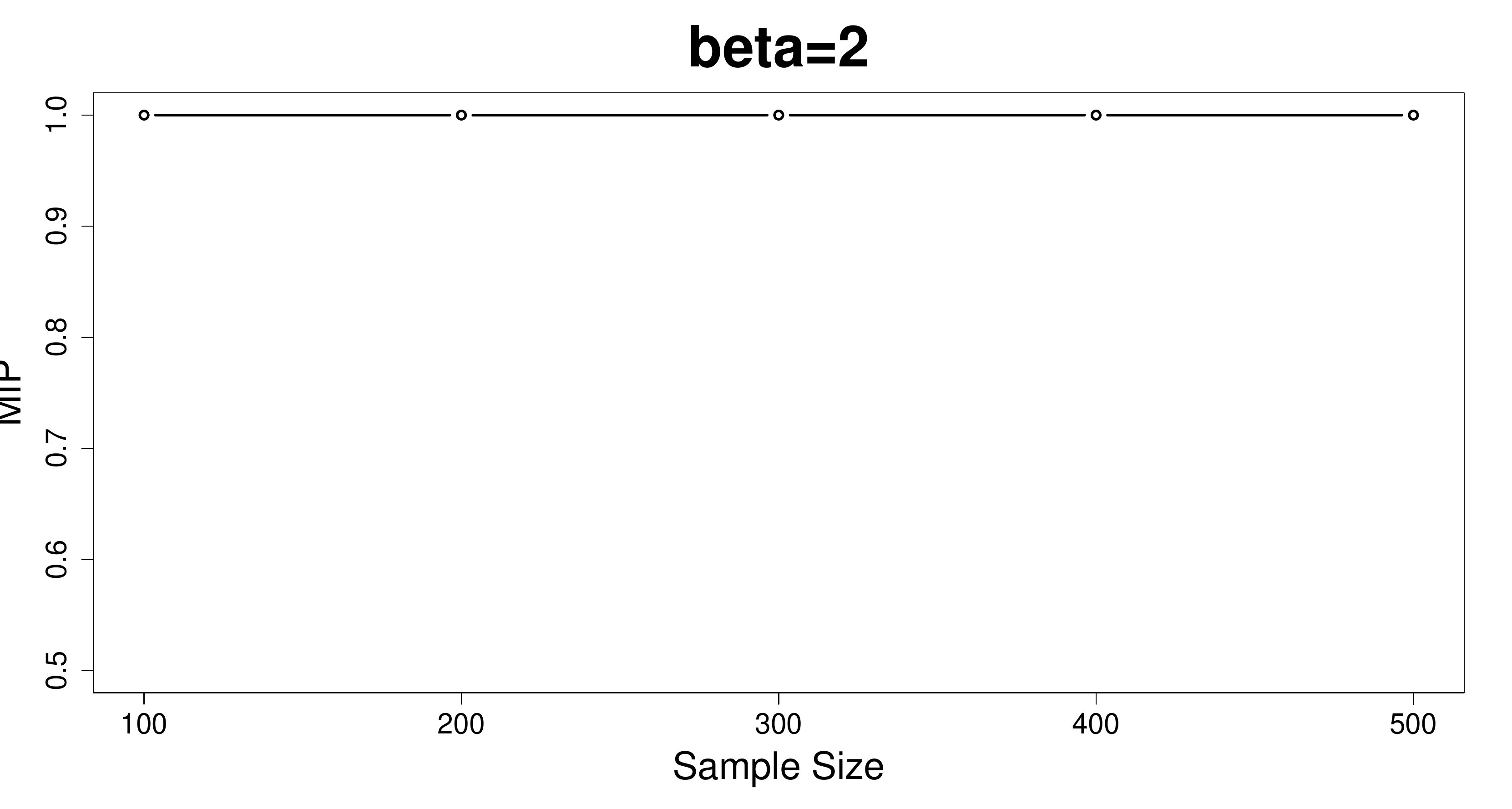}} \quad \mbox{\includegraphics[height=2.0in, width=0.45\textwidth]{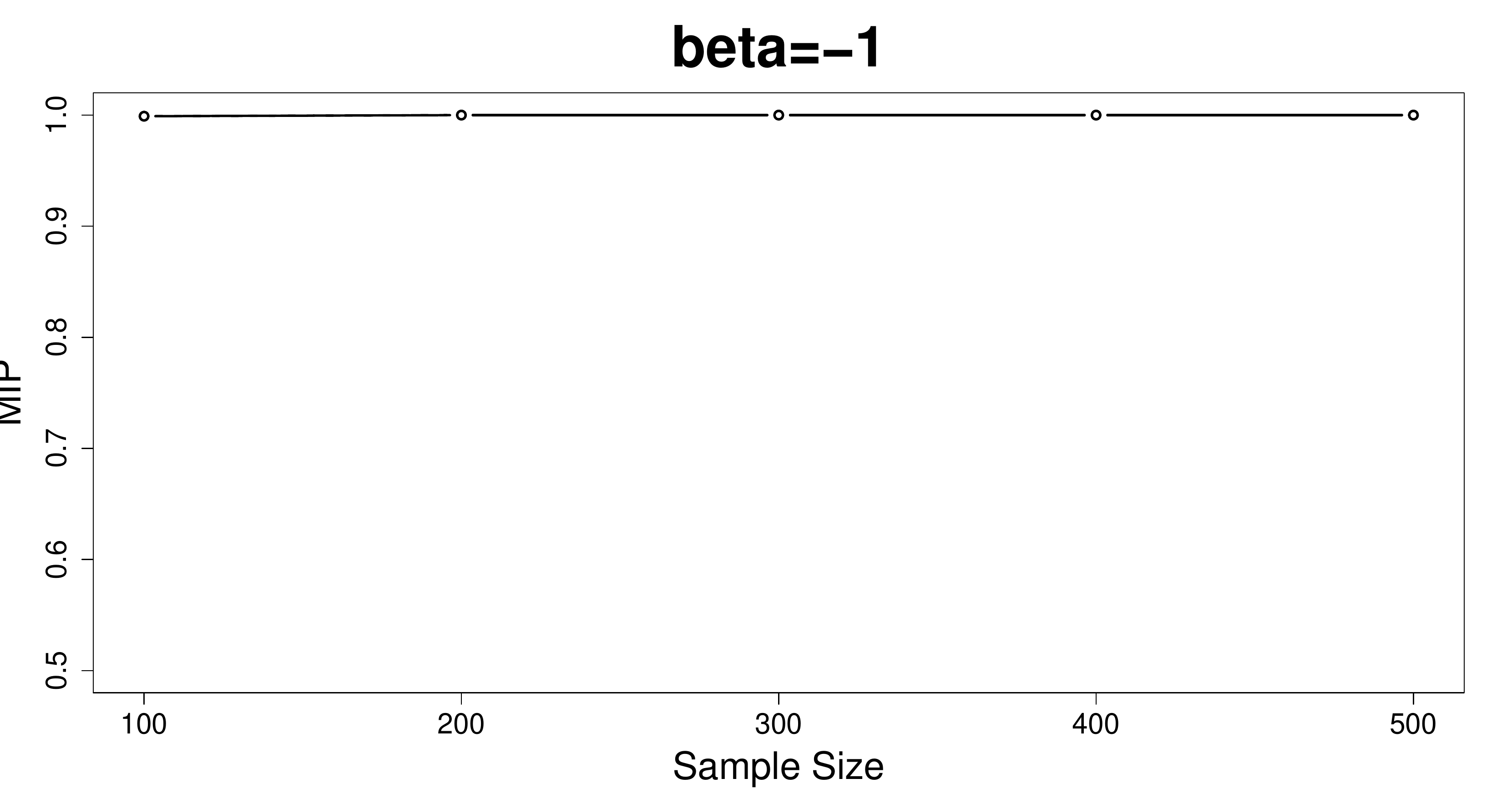}}\\
		\mbox{\includegraphics[height=2.0in, width=0.45\textwidth]{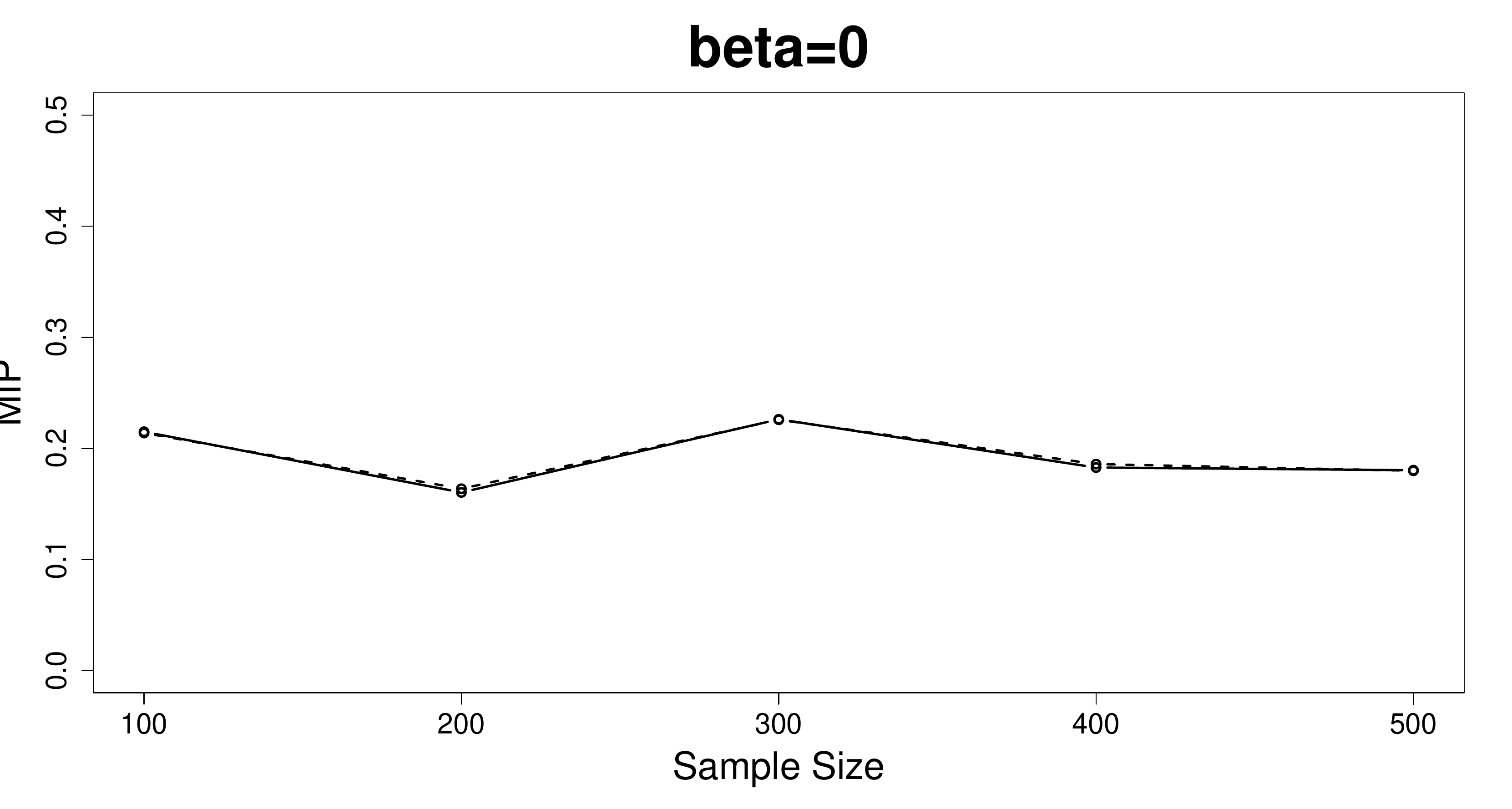}} \quad 	\mbox{\includegraphics[height=2.0in, width=0.45\textwidth]{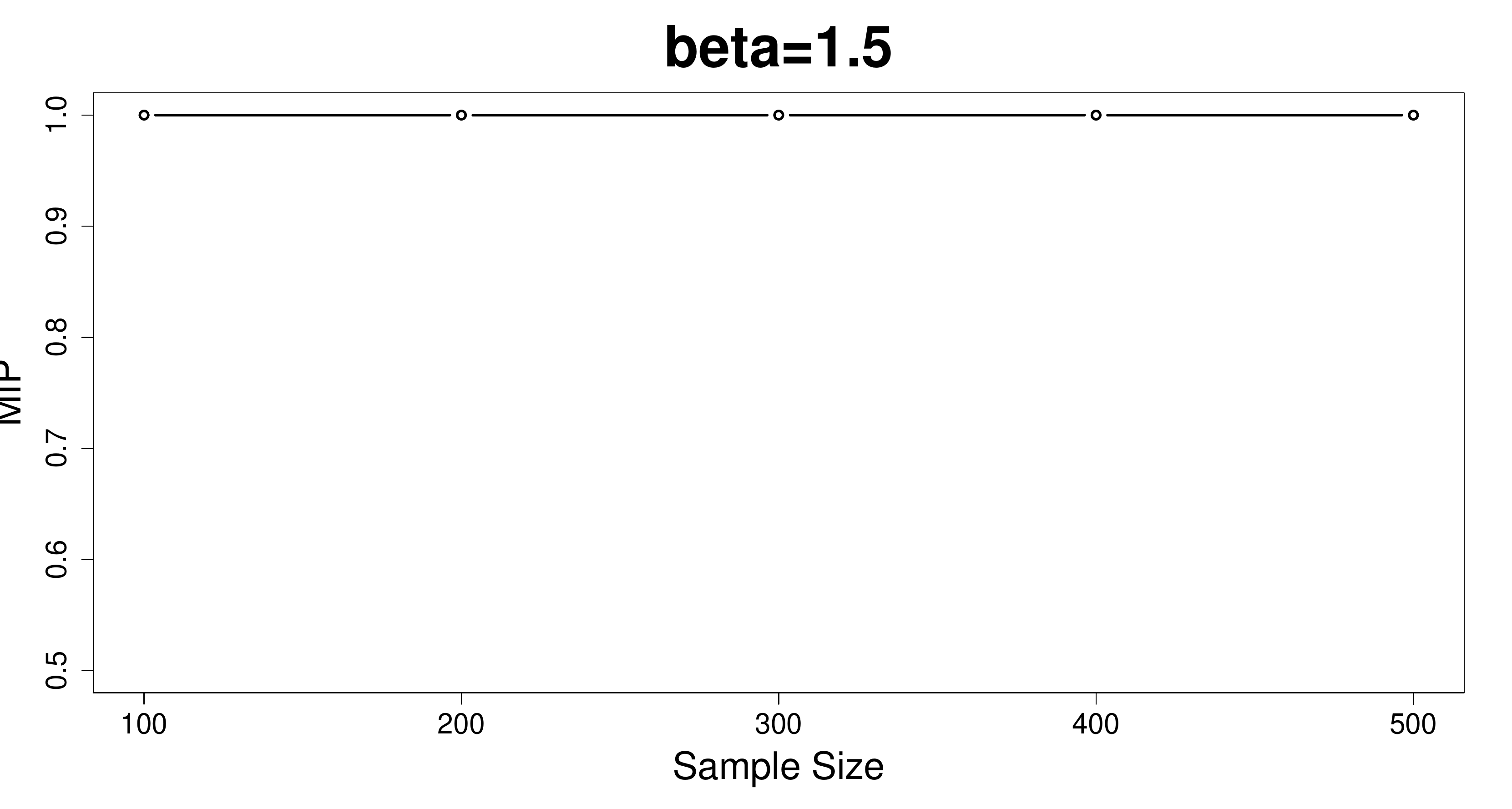}}\\
			\mbox{\includegraphics[height=2.0in, width=0.45\textwidth]{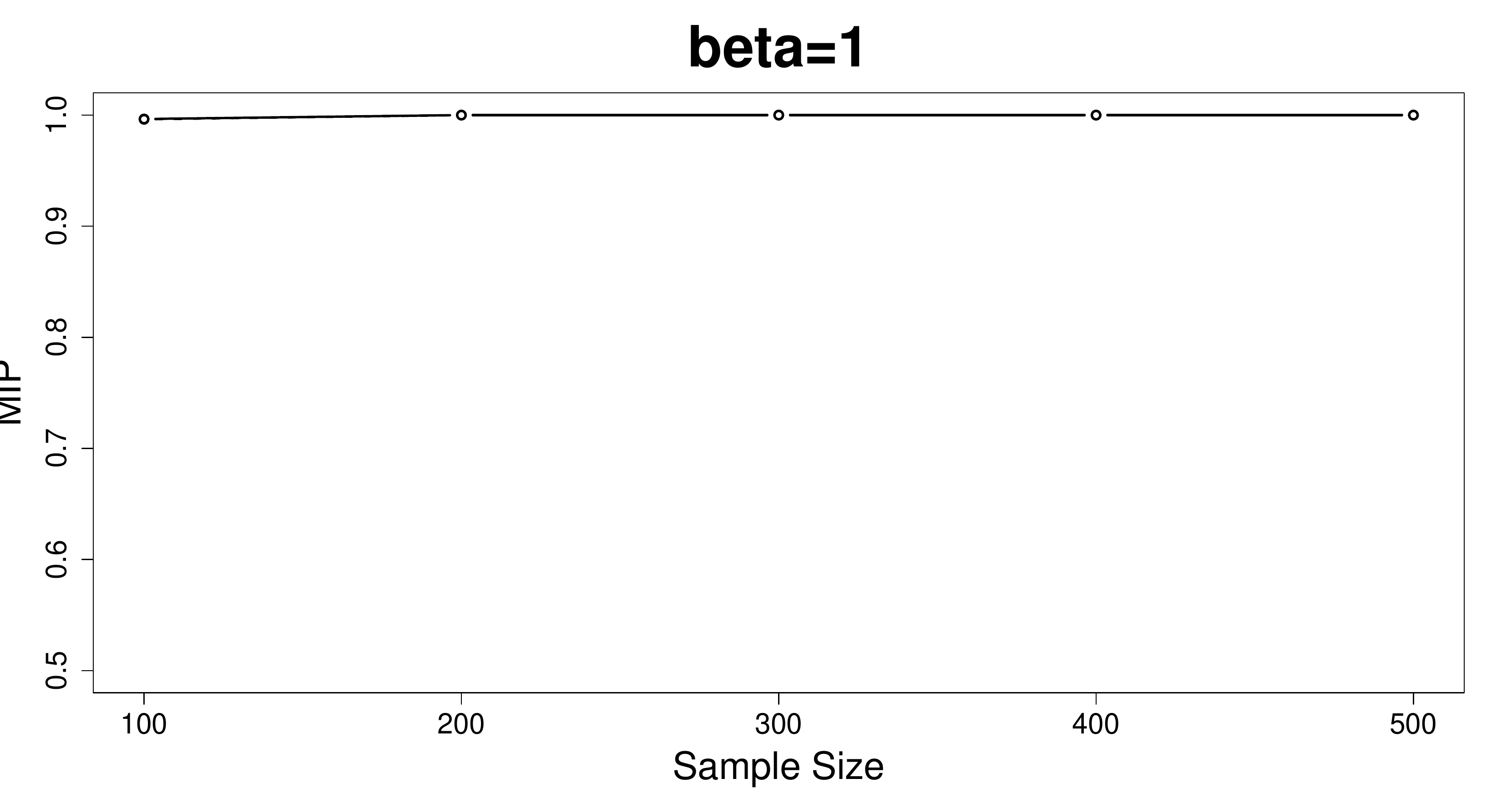}} \quad \mbox{\includegraphics[height=2.0in, width=0.45\textwidth]{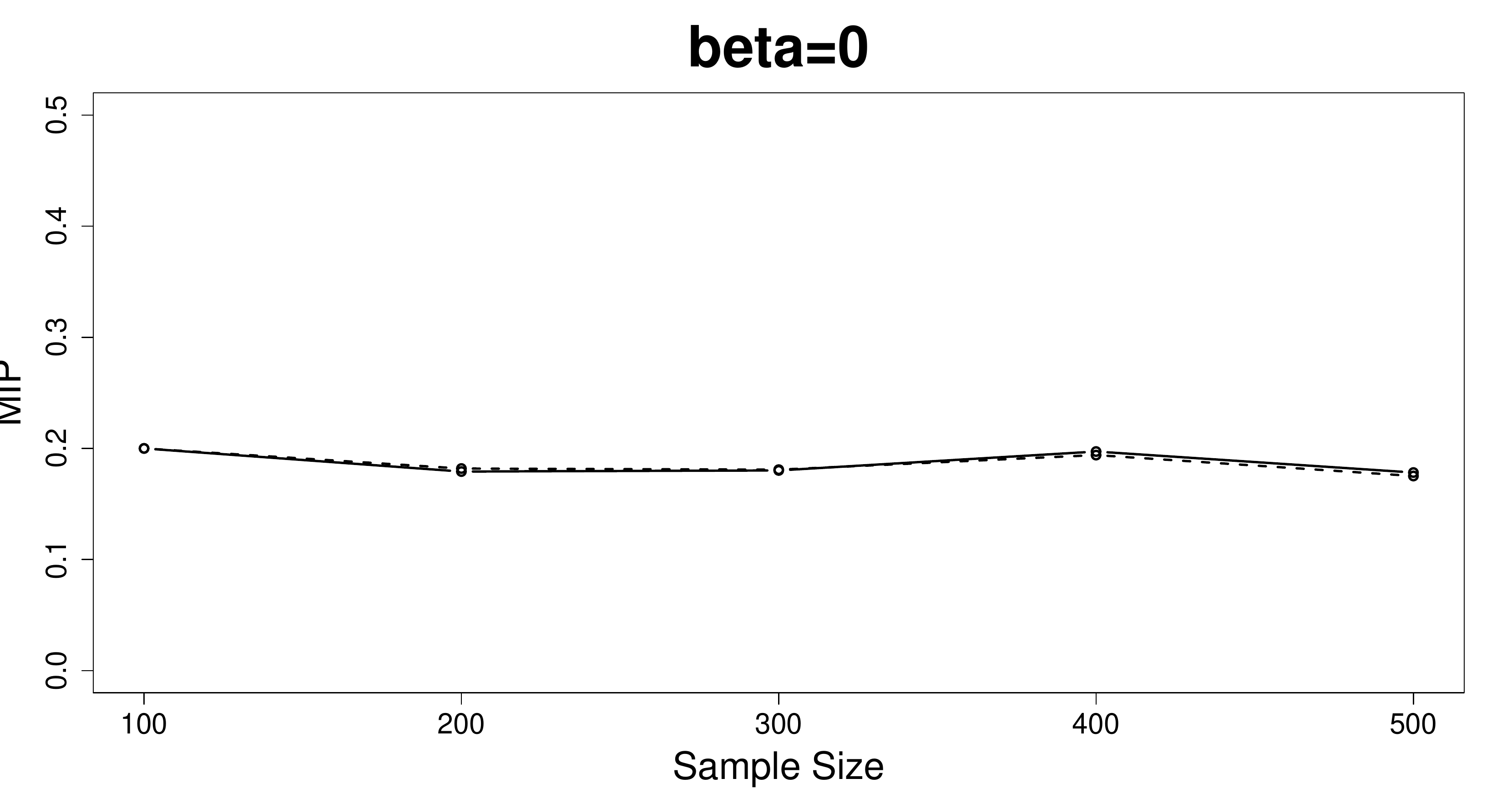}}\\
		\mbox{\includegraphics[height=2.0in, width=0.45\textwidth]{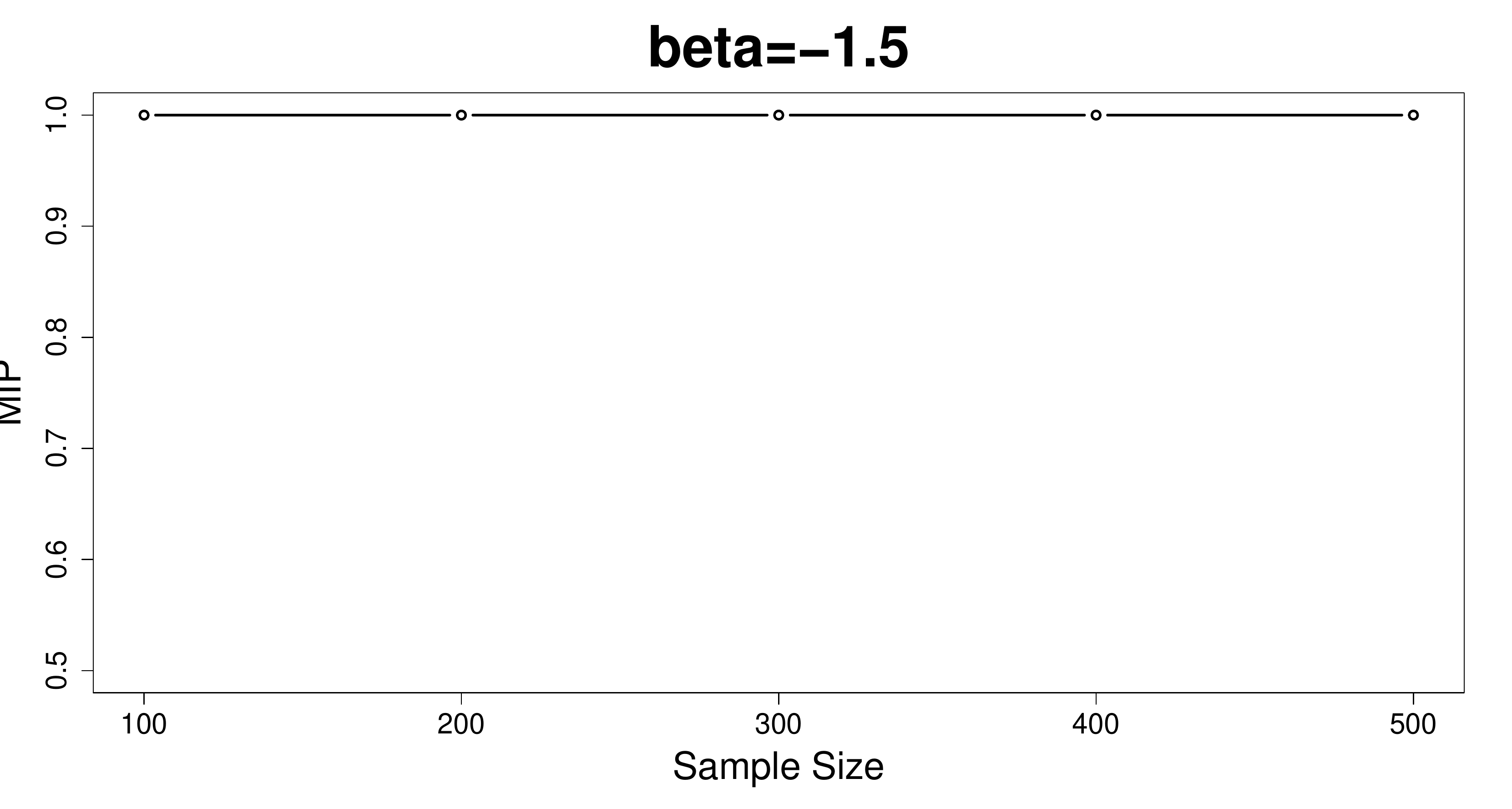}} \quad 	\mbox{\includegraphics[height=2.0in, width=0.45\textwidth]{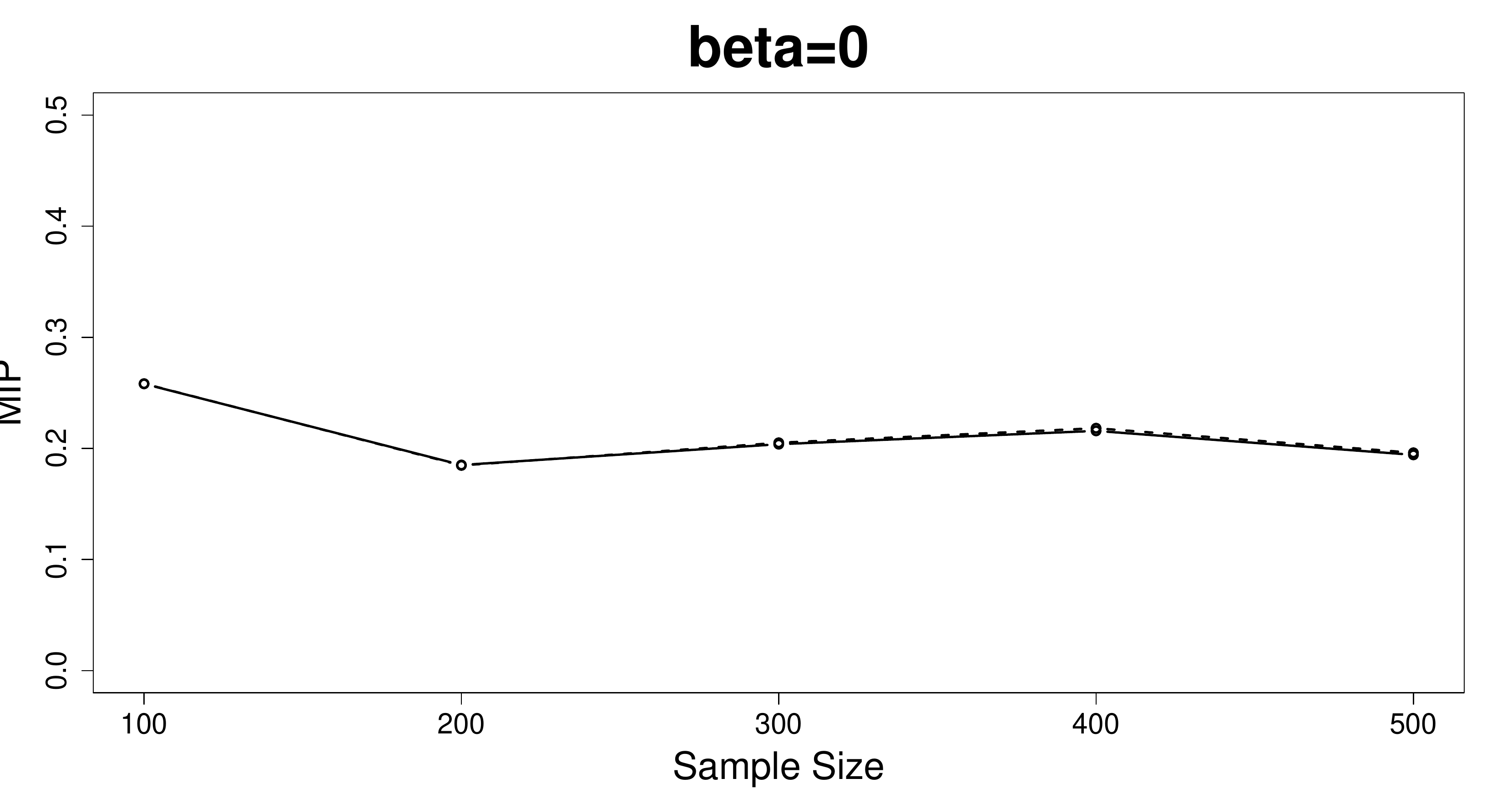}}
			\caption{\small{Marginal Inclusion Probabilities (MIP): Truth generated from Gaussian residual. Solid lines - Semi-parametric Linear Model, dashed lines - Non-parametric Linear Model.}}
\end{figure}

\begin{figure}
		\mbox{\includegraphics[height=3.0 in, width=1 \textwidth]{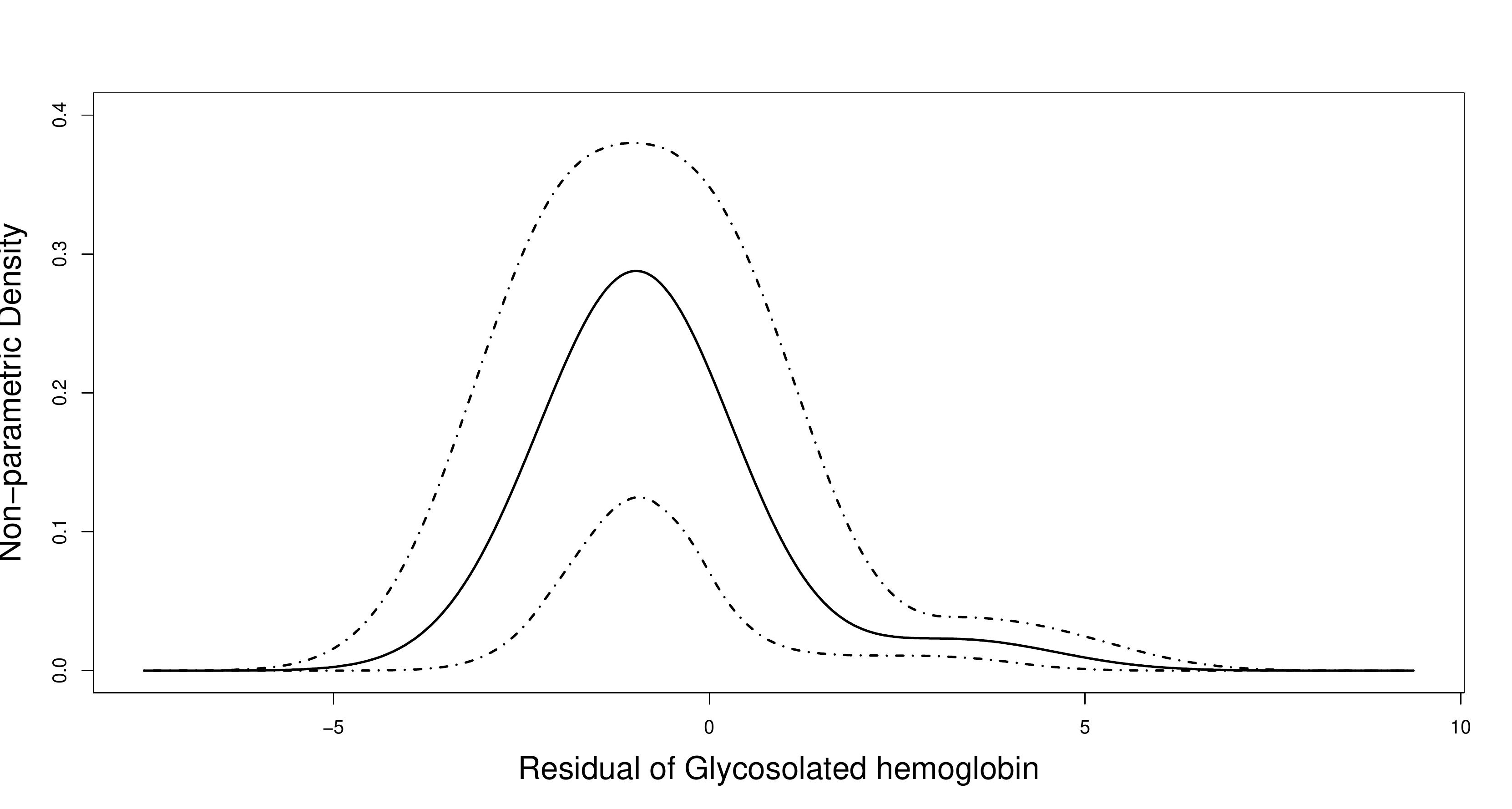}} 		
		\caption{\small{Residual plots for Diabetes study for Semi-parametric Linear Model}}
\end{figure}


\begin{thebibliography}{}

\bibitem{ } Bartlett, M. (1957), ``A comment on D. V. Lindley's statistical paradox", Biometrika, 44, 533-534. 
\bibitem{ } Berger, J. O. and Pericchi, L. R. (1996), ``The intrinsic Bayes factor for model selection and prediction", J. Amer. Statist. Assoc., 91, 109 - 122.
\bibitem{ } Berger, J. O. and Pericchi, L. (2001), ``Objective Bayesian methods for model selection: Introduction
and comparison", Model Selection, vol. 38 of IMS Lecture Notes - Monograph Series, 135 - 193. Institute of Mathematical Statistics. 
\bibitem{ } Brunham, L.R., Kruit, J.K., Pape, T.D., Timmins, J.M., Reuwer, A.Q., Vasanji, Z., Marsh, B.J., Rodrigues, B., Johnson, J.D., Parks, J.S., Verchere, C.B., and Hayden, M.R. (2007), ``$\beta$-cell ABCA1 influences insulin secretion, glucose homeostasis and response to thiazolidinedione treatment", Nature Medicine, 13, 340 - 347. 
\bibitem{ } Casella, G. and Moreno E. (2006), ``Objective Bayesian variable selection", J. Amer. Statist. Assoc., 101, 157 - 167. 
\bibitem{ } Casella, G., Gir$\acute{o}$n, F. J., Mart$\acute{i}$nez, M. L. and Moreno, E. (2009), ``Consistency of
Bayesian procedures for variable selection", Ann. Statist. 37 1207 - 1228. 
\bibitem{ } Chung, Y., and Dunson, D.B. (2009), ``Nonparametric Bayes Conditional Distribution Modeling With Variable Selection", J. Amer. Statist. Assoc., 104(488), 1646-1660. 
\bibitem{ } Dunson, D.B., and Herring, A.H. (2005), ``Bayesian model selection and averaging in additive and proportional hazards models", Lifetime Data Analysis, 11(2), 213-232. 
\bibitem{ } Epstein, M., and Sowers, J.R. (1992), ``Diabetes mellitus and hypertension", Hypertension, 19, 403-418. 
\bibitem{ } Fern$\acute{a}$ndez, C., E. Ley, and M. F. J. Steel (2001), `'Benchmark priors for
Bayesian model averaging", J. Econometrics, 100(2), 381 - 427. 
\bibitem{ } George, E. I. and McCulloch, R. E. (1993), ``Variable Selection Via Gibbs Sampling", J. Amer. Statist. Assoc., 88(423), 881-89. 
\bibitem{ } George, E. I. and McCulloch, R. E. (1997), ``Approaches for Bayesian Variable Selection", Statist. Sinica, 7(2), 339-74.
\bibitem{ } Ghosal, S., Lember, J., and van der Vaart, A. (2008), ``Nonparametric Bayesian model selection and averaging",  Electronic J. Stat., 2, 63-89. 
\bibitem{ } Green, P.J.  (1995), ``Reversible jump Markov chain Monte Carlo computation and Bayesian model determination", Biometrika, 82 (4), 711-732. 
\bibitem{ } Guo, R. and Speckman, P. (2009), ``Bayes factor consistency in linear models", In the 2009 International Workshop on Objective Bayes Methodology, Philadelphia, 2009.
\bibitem{ } Jeffreys, H. (1961), ``Theory of Probability", Oxford Univ. Press. 
\bibitem{ } Jiang, W. (2007), `` Bayesian Variable Selection for high dimensional generalized linear models: convergence rates of the fitted densities'', Ann. of Statist., 35(4), 1487 - 1511.
\bibitem{ } Kass, R. E.,and Raftery, A.E.(1995), ``Bayes Factors", J. Amer. Statist. Assoc., 90, 773 - 795.
\bibitem{ } Kim, S., Tadesse, M.G., and Vannucci, M. (2006), ``Variable selection in clustering via Dirichlet process mixture models", Biometrika, 93(4), 877-893. 
\bibitem{ } Kuo, L. and Mallick, B. (1997), ``Semiparametric inference for the accelerated failure time model", Can. J. Stat., 25, 457-472.
\bibitem{ } Kyung, M., Gill, J., and Casella, G. (2009), ``Characterizing the variance improvement in linear
Dirichlet random effects models", Statistics and Probability Letters, 79, 2343-2350. 
\bibitem{ } Liang, F., Paulo, R., Molina, G., Clyde, M.A., and Berger, J.O. (2008), ``Mixtures of g-priors for Bayesian Variable Selection.", J. Amer. Statist. Assoc., 103(481), 410-423.
\bibitem{ } Marin, J.M. and Robert, C. P. (2007), ``Bayesian Core: A Practical Approach to Computational Bayesian Statistics", Springer-Verlag Inc. 
\bibitem{ } Meyer, M. C. and Laud, P. W. (2002), ``Predictive variable selection in generalized linear models", J. Amer. Statist. Assoc., 97, 859 - 871.
\bibitem{ } Mokdad, A.H., Bowman, B.A., Ford, E.S., Vinicor, F., Marks, J.S., Koplan, J.P. (2001), ``The continuing epidemics of obesity and diabetes in the United States", J. Amer. Med. Assoc., 286, 1195-1200. 
\bibitem{ } Moreno, E., Bertolino, F. and Racugno, W. (1998), ``An intrinsic limiting procedure for model selection and hypothesis testing", J. Amer. Statist. Assoc., 93, 1451 - 1460.
\bibitem{ } Moreno, E., Gir$\acute{o}$n, F.J., and Casella, G. (2010), ``Consistency of objective Bayes factors as the model dimension grows", Ann. Statist., 38(4), 1937 - 1952. 
\bibitem{ } Mostofi, A.G., and Behboodian, J. (2007), ``On model selection in Bayesian regression", Metrika, 66(3), 259-268. 
\bibitem{ } O'Hara, R.B. and Sillanp\"a\"a, M.J. (2009), ``Review of Bayesian variable selection methods: What, how and which", Bayesian Analysis, 4, 85 - 118.
\bibitem{ } Park, T., and Casella, G. (2008), ``The Bayesian Lasso", J. Amer. Statist. Assoc., 103(482), 681-686. 
\bibitem{ } Raftery, A. E. and Richardson, S. (1993), ``Model selection for generalized linear
models via GLIB, with application to epidemiology", Bayesian Biostatistics, Berry, D. A. and Stangl, D. K., editors. Marcel Dekker, New York.
\bibitem{ } Ritter, C., and Tanner, M.A. (1992), ``Facilitating the Gibbs sampler: the Gibbs stopper and the griddy-Gibbs sampler", J. Amer. Statist. Assoc., 87(419), 861-868.
\bibitem{ } Schmidt, M.I., Duncan, B.B., Canani L.H., Karohl, C., and Chambless L. (1992), ``Association of waist-hip ratio with diabetes mellitus. Strength and possible modifiers", Diabetes Care. 15(7), 912-4.
\bibitem{ } Smith, M., and Kohn, R. (1996), ``Nonparametric regression using Bayesian variable selection", J. Econometrics, 75(2), 317-343. 
\bibitem{ } Shao, J. (1997), ``An asymptotic theory for linear model selection", Statist. Sinica, 7, 221 - 264. 
\bibitem{ } Strawderman, W. E. (1971), ``Proper Bayes minimax estimators of the multivariate normal mean", Ann. Math. Statist., 42, 385-388. 
\bibitem{ } Tibshirani, R. (1996), ``Regression shrinkage and selection via the lasso", J. Royal. Statist. Soc., Series B., 58(1), 267-288. 
\bibitem{ } Walker, S. (2007), ``Sampling the dirichlet mixture model with slices", Comm. Statist. Sim. Comput., 36, 45 - 54. 
\bibitem{ } Yau C., Papaspiliopoulos, O., Roberts, G. and Holmes, C. (2011), `'Bayesian non-parametric hidden Markov models with applications in genomics",  J. Royal Stat. Soc., Series B, 73(Part 1), 33 - 57. 
\bibitem{ } Yi, N., and S. Xu. (2008), ``Bayesian LASSO for quantitative trait loci mapping", Genetics, 179, 1045-1055. 
\bibitem{ } Zellner, A. and Siow, A. (1980), ``Posterior odds ratios for selected regression hypotheses", Bayesian Statistics: Proceedings of the First International Meeting, Valencia, University of Valencia Press, 585-603.
\bibitem{ } Zellner, A. (1986), ``On assessing prior distributions and Bayesian regression analysis with g-prior
distributions", Bayesian Inference and Decision Techniques: Essays in Honor of Bruno de
Finetti, 233-243. 
\bibitem{ } Zou, H., and Hastie, T. (2005), ``Regularization and variable selection via the elastic net", 
J. Royal. Statist. Soc, Series B., 67(2), 301 - 320. 
\end{thebibliography}
\end{document}